\documentclass[11pt]{amsart}
  \usepackage{amsmath, amsfonts, amssymb,amsthm}

\newtheorem{theorem}{Theorem}
\newtheorem{lemma}{Lemma}
\newtheorem{cor}{Corollary}

\theoremstyle{definition}

\def\min{\mathop{\mathrm{min}}}

\newcommand{\C}{\mathbb{C}}

\newcommand{\N}{\mathbb{N}}

\newcommand{\PP}{\mathbb{P}}

\newcommand{\Z}{\mathbb{Z}}

\newcommand{\ax}{\rightarrow }

\newcommand{\ab}{{\bf a}}

\begin{document}
\title{{Difference analogue of second main theorems for meromorphic mapping into algebraic variety}}

\author{Pei-Chu Hu$^1$ and Nguyen Van Thin$^{1,2}$}
\address{Shandong University$^{1}$, Department of Mathematics, Jinan 250100, Shandong Province, P. R. China}
\address{Thai Nguyen University of Education$^{2}$, Department of Mathematics, Luong Ngoc Quyen  street, Thai Nguyen city, Thai Nguyen, Viet Nam.}
\email{pchu@sdu.edu.cn}
\email{thinmath@gmail.com}


\begin{abstract}
In this paper, we prove some difference analogue of second main theorems of meromorphic mapping from $\mathbb C^m$ into an algebraic variety $V$
  intersecting a finite set of fixed hypersurfaces in subgeneral position. As an application, we prove a result on
 algebraically degenerate of holomorphic curves on $\mathcal P^1_{c}$ intersecting hypersurfaces and difference analogue of Picard's theorem on holomorphic curves. Furthermore, we obtain a second main theorem of meromorphic mappings intersecting hypersurfaces in $N$-subgeneral position for Veronese embedding in $\mathbb P^n(\mathbb C)$ and a uniqueness theorem sharing hypersurfaces.
\end{abstract}
\let\thefootnote\relax\footnote
{The research results are sponsored by China/Shandong University International Postdoctoral Exchange Program, NSFC of Shandong (No. ZR2018MA014), PCSIRT (No. IRT1264) and The Fundamental Research Funds of Shandong University(No. 2017JC019).}
\baselineskip=16truept
\maketitle
\pagestyle{myheadings}
\markboth{}{}

\section{ Introduction and main results}
\def\theequation{1.\arabic{equation}}
\setcounter{equation}{0}
Recently, the second main theorem of Nevanlinna have been studied actively for difference operators. For example, R. Halburd and R. Korhonen \cite{RRK1, RRK2} in 2006 built the second main theorem for difference operators of meromorphic functions on $\mathbb{C}$. In 2014, R. Halburd, R. Korhonen and K. Tohge \cite{RRK} proved the difference analogue  of second main theorem of holomorphic curves from $\mathbb C$  into $\mathbb P^n(\mathbb C)$ intersecting a finite set of fixed hyperplanes in general position.

In 2009, M. Ru \cite{Ru3} proved the second main theorem of holomorphic curves into an algebraic variety. In 2017, S. D. Quang \cite{Q} extended the result of M. Ru  \cite{Ru3} to hypersurfaces in subgeneral position. Our goal is to establish the difference analogue  of second main theorem of meromorphic mappings from $\mathbb C^m$  into an algebraic variety $V$ intersecting a finite set of fixed hypersurfaces in subgeneral position depending on a step number of difference. In particular, for the Veronese embedding in $\mathbb P^n(\mathbb C),$ our second main theorem and difference analogue  of  Picard's theorem recover the results of Cao-Korhonen \cite{TC} and Halburd-Korhonen-Tohge \cite{RRK}, respectively. By a way, we also obtain uniqueness theorems of meromorphic mappings which improve the result of Dulock-Ru \cite{DR}.

To introduce our results clearly, it is necessary to introduce some notations. Take $z=(z_1, \dots, z_m)\in \mathbb C^m$ and write the standard norm $||z||=(\sum_{j=1}^{m}|z_j|^2)^{1/2}$. For $r>0,$ define the ball $B(r)=\{z\in \mathbb C^m: ||z||<r\}$ and the sphere $S(r)=\{z\in \mathbb C^m: ||z||=r\}$. As usual, define the differential operator $d^c=\dfrac{1}{4\pi i}(\partial-\overline \partial)$ and two differential forms $v(z)=(dd^{c}||z||^2)^{m-1}$, $\sigma(z)=d^c\log ||z||^2\wedge (dd^c||z||^2)^{m-1}.$

Take $z_0\in \mathbb C^m$ and  $a\in \mathbb P^1(\mathbb C)$. We expand a nonzero entire function $h$ on $\mathbb C^m$ into a Taylor series $h(z)=\sum_{n=0}^{\infty}P_n(z-z_0)$ at $z_0$,
 where $P_n(z)$ is either identically zero or a homogeneous polynomial of $z$ with degree $n.$ The number $\nu_h(z_0)=\min \{n: P_n \ne 0\}$ is said to be the zero multiplicity of $h$ at $z_0.$ Set ${\rm supp } \nu_h:=\overline{ \{z\in \mathbb C^m: \nu_h(z)\ne 0\}},$ which is a purely $(m-1)$-dimensional analytic subset or empty set. Further, if $h$ is a  nonzero meromorphic function on $\mathbb C^m$, we can choose two holomorphic functions $h_0$ and $h_1$ on a neighborhood $U$ of $z_0$ such that $h=\dfrac{h_0}{h_1}$ on $U$ and $\dim h_0^{-1}(0)\cap h_1^{-1}(0)\le m-2$. Define  the $a$-valued multiplicity of $h$ at $z_0$ by $\nu_{h-a}(z_0)=\nu_{h=a}(z_0):=\nu_{h_0-ah_1}(z_0)$, where particularly $\nu_{h}=\nu_{h=0}=\nu_{h_0}, \nu_{h=\infty}=\nu_{h_1},$ which are independent of the choices of $h_0$ and $h_1.$

In Nevanlinna's theory, a multiplicity function (or divisor) $\nu$ on $\mathbb C^m$ usually be associated with the following truncation functions
$$ \nu^{M}(z)=\min \{\nu(z), M\},\;
\nu_{ \le k}^{M}(z)=
\begin{cases}
  & 0\hspace{1cm}\text{if}\; \nu(z)>k\\
&\nu^{M}(z)\; \text{if}\; \nu (z)\le k
\end{cases}
$$
and
$$ \nu_{ \ge k}^{M}(z)=
\begin{cases}
  & 0\hspace{1cm}\text{if}\; \nu(z)<k\\
&\nu^{M}(z)\; \text{if}\; \nu (z)\ge k
\end{cases},$$
where $k, M$ are positive integers or $+\infty.$ Further, the function $\nu$ defines a counting function as follows:
$$ n(t)=
\begin{cases}
&\int_{{\rm supp} \nu\cap B(t)}\nu v \hspace{0.6cm}\text{if}\; m\ge 2\\
&\sum_{|z|\le t}\nu(z)\hspace{1cm} \text{if}\; m=1
\end{cases}.$$
Similarly, the truncation functions $\nu^{M}$, $\nu_{ \le k}^{M}$ and $\nu_{ \ge k}^{M}$ can define counting functions $n^{M}(t), n_{\le k}^{M}(t)$ and $n_{\ge k}^{M}(t)$, respectively. Moreover, the function  $\nu$ also defines a valence function
$$ N(r, \nu)=\int\limits_{1}^{r}\dfrac{n(t)}{t^{2m-1}}dt$$
for $r>1$, so that the valence functions $N(r, \nu^{M}), N(r, \nu_{\le k}^{M})$, $N(r, \nu_{\ge k}^{M})$ are associated with the truncation functions $\nu^{M}$, $\nu_{ \le k}^{M}$ and $\nu_{ \ge k}^{M}$, respectively, which also be denoted by $N^{M}(r, \nu),$
 $N_{\le k}^{M}(r, \nu)$ and $N_{\ge k}^{M}(r, \nu),$ respectively. In particular, for a meromorphic function $h$ on $\mathbb C^m,$ we write
$$N_{h}(r)=N(r, \nu_{h}), N_{h}^{M}(r)=N(r, \nu_{h}^{M}), $$
$$N_{h, \le k}^{M}(r)=N(r, \nu_{h, \le k}^{M}), N_{h, \ge k}^{M}(r)=N(r, \nu_{h, \ge k}^{M}).$$

Let $f:\mathbb C^m\longrightarrow \mathbb P^n(\mathbb C)$ be a nonconstant meromorphic mapping. We can choose holomorphic functions $f_0, f_1, \dots, f_n$ on $\mathbb C^m$ such that $I_f:=\{z\in \mathbb C^m: f_0(z)=\dots=f_n(z)=0\}$ is of dimension at most $m-2$ and $f=(f_0:\dots:f_n).$ Usually, $\tilde{f}=(f_0,\dots,f_n):\mathbb{C}^m\longrightarrow \mathbb{C}^{n+1}$ is called a reduced representation of $f.$ For $r>1$, the characteristic function of $f$ can be given by
$$ T_f(r)=\int\limits_{S(r)}\log ||\tilde{f}(z)||\sigma(z)-\int\limits_{S(1)}\log ||\tilde{f}(z)||\sigma(z),$$
 where $||\tilde{f}(z)||=(\sum_{j=0}^{n}|f_j(z)|^{2})^{1/2}.$  Further, we denote the {\it hyper-order} and {\it order} of $f$ respectively by
 $$\varsigma(f)=\limsup_{r\to\infty}\dfrac{\log \log T_f(r)}{\log r},\quad \sigma(f)=\limsup_{r\to\infty}\dfrac{\log T_f(r)}{\log r}.$$
Moreover, we use $S(r, f)$ (resp., $S_1(r,f)$)  to denote a quantity satisfying the following property
\begin{equation*}
    S(r,f)=o(T_f(r))
\end{equation*}
outside a set of finite logarithmic (or Lebesgue) measure (resp., $S_1(r,f)=o(T_f(r))$ on a set of logarithmic density $1$).

Take ${\bf a}=(a_0,\dots,a_n)\in\mathbb{C}^{n+1}$ and define a linear form
\begin{equation*}
   L(x_0,\dots, x_n)=  L_{{\bf a}}(x_0,\dots, x_n)=\sum_{j=0}^{n}a_jx_j
\end{equation*}
associated with a hyperplane
\begin{equation*}
H=H_{\bf a}=\{x=(x_0:\dots:x_n)\in \mathbb P^n(\mathbb C)\ |\ L_{{\bf a}}(x)=L_{{\bf a}}(x_0,\dots, x_n)=0\}.
\end{equation*}
If $f(\mathbb C^m)\not\subset H,$ the proximity function of $f$ with respective to $H$ can be obtained by
$$ m_{f}(r, H) =\int\limits_{S(r)}\log \dfrac{||\tilde{f}(z)||\cdot||L_{\bf a}||}{|L_{\bf a}(\tilde{f}(z))|}\sigma(z)-\int\limits_{S(1)}\log \dfrac{||\tilde{f}(z)||\cdot||L_{\bf a}||}{|L_{\bf a}(\tilde{f}(z))|}\sigma(z)$$
for $r>1$, where $||L_{\bf a}||=(\sum_{j=0}^{n}|a_j|^2)^{1/2}$, in which, up to an additive constant, we may choose $||\tilde{f}(z)||=\max\{|f_0(z)|, \dots, |f_n(z)|\}$ and $||L_{\bf a}||=\sum_{j=0}^{n}|a_j|$ since norms on $\mathbb C^m$ are equivalent.

Generally,  if $Q(x_0,\dots,x_n)\in \mathbb{C}[x_0,...,x_n]$ is a homogeneous polynomial of degree $d$, then a hypersurface
\begin{equation*}
    D=\{x=(x_0:\dots:x_n)\in \mathbb P^n(\mathbb C)\ |\ Q(x)=Q(x_0,\dots, x_n)=0\}
\end{equation*}
of degree $d$ in $\mathbb P^n(\C)$ is associated with $Q$. If $f(\mathbb C^m)\not\subset D,$ i.e., $Q(\tilde{f}) \not\equiv 0,$ the proximity function $m_f(r, D)$ of $f$ can be given by
$$m_{f}(r,D)=\int\limits_{S(r)} \log \dfrac {||\tilde{f}(z)||^{d}||Q||}{|Q(\tilde{f}(z))|}\sigma(z)-\int\limits_{S(1)} \log \dfrac {||\tilde{f}(z)||^{d}||Q||}{|Q(\tilde{f}(z))|}\sigma(z)$$
for $r>1$, where $||Q||$ is the total of the absolute values of the coefficients of $Q.$
Further, let $\nu_{Q(\tilde{f})}$ be the zero multiplicity function of $Q(\tilde{f})$. Then the following valence functions
$$N_{f}(r, D)=N(r, \nu_{Q(\tilde{f})}), N_{f}^{M}(r, D)=N(r, \nu_{Q(\tilde{f})}^{M}), $$
$$N_{f, \le k}^{M}(r, D)=N(r, \nu_{Q(\tilde{f}), \le k}^{M}), N_{f, \ge k}^{M}(r, D)=N(r, \nu_{Q(\tilde{f}), \ge k}^{M})$$
are defined well. We also use the symbol $\overline E_{M)}(D, f)$ to denote ${\rm supp}\nu_{Q(\tilde{f}), \le M}^{M}$, that is, $\overline E_{M)}(D, f)$ is the set of zeros of $Q(\tilde{f})$ with multiplicity at most $M,$ in which each zero is counted only one time.

Our work is based on a decomposition of the $a$-valued multiplicity $\nu_{h-a}$ of a meromorphic function $h$ on $\mathbb C^m$
\begin{equation*}
\nu_{h-a}=    \nu_{h-a}^{[M,c]}+\widetilde{\nu}_{h-a}^{[M,c]}
\end{equation*}
defined by R. Korhonen, N. Li and K. Tohge in \cite{RRK3}, in which $M$ is a positive integer, $c\in\mathbb{C}^m\setminus\{0\}$, and
\begin{equation*}
  \nu_{h-a}^{[M,c]}(z)=\left\{
                      \begin{array}{ll}
                        \nu_{h-a}(z), & \hbox{ if $\nu_{h_{jc}-a}(z)\geq \nu_{h-a}(z)$ for $j=1,...,M$;} \\
                        0, & \hbox{otherwise,}
                      \end{array}
                    \right.
\end{equation*}
where $h_{jc}(z)=h(z+jc)$. The points on ${\rm supp}\nu_{h-a}^{[M,c]}$ (resp., ${\rm supp}\widetilde{\nu}_{h-a}^{[M,c]}$) are called to be  $M-$ successive and $c-$separated (resp., $M$-aperiodic of pace $c$), which defines the valence function $N_{h-a}^{[M,c]}(r)$ (resp., $\overset{\sim} N_{h-a}^{[M, c]}(r)$). It is obvious that the number $M$ serves as the step of difference.

Thus the zero multiplicity function $\nu_{Q(\tilde{f})}$ of $Q(\tilde{f})$ defined by the hypersurface $D$ in $\mathbb P^n(\C)$ of degree $d$ corresponds Korhonen-Li-Tohge's decomposition
\begin{equation*}
    \nu_{Q(\tilde{f})}=\nu_{Q(\tilde{f})}^{[M,c]}+\widetilde{\nu}_{Q(\tilde{f})}^{[M,c]}
\end{equation*}
associated with valence functions $N_{f}^{[M, c]}(r, D)$ and $\overset{\sim} N_{f}^{[M, c]}(r, D)$ of $\nu_{Q(\tilde{f})}^{[M,c]}$ and $\widetilde{\nu}_{Q(\tilde{f})}^{[M,c]}$, respectively, and hence
$$\overset{\sim} N_{f}^{[M, c]}(r, D)=\overset{\sim} N_{Q(\tilde{f})}^{[M, c]}(r).$$
Similarly, we also consider the following decomposition
\begin{equation*}
    \nu_{Q(\tilde{f})}=\dot{\nu}_{Q(\tilde{f})}^{[M,p]}+\hat{\nu}_{Q(\tilde{f})}^{[M,p]}
\end{equation*}
in which $p=(p_1, \dots, p_m)\in\mathbb{C}^m\setminus\{0,1\}$, and
\begin{equation*}
 \dot{ \nu}_{Q(\tilde{f})}^{[M,p]}(z)=\left\{
                      \begin{array}{ll}
                        \nu_{Q(\tilde{f})}(z), & \hbox{ if $\nu_{Q(\tilde{f})_{p^{j}}}(z)\geq \nu_{Q(\tilde{f})}(z)$ for $j=1,\dots,M$;} \\
                        0, & \hbox{otherwise,}
                      \end{array}
                    \right.
\end{equation*}
where $Q(\tilde{f})_{p^{j}}(z)=Q(\tilde{f})(p^{j}\cdot z)=Q(\tilde{f})(p_1^{j}z_1, \dots, p_m^{j}z_m)$, and denote valence functions $\dot{N}_{f}^{[M, p]}(r, D)$ and $\hat{N}_{f}^{[M, p]}(r, D)$ of $\dot{\nu}_{Q(\tilde{f})}^{[M,p]}$ and $\hat{\nu}_{Q(\tilde{f})}^{[M,p]}$, respectively. Now we can state first result as follows:

\begin{theorem}\label{th1}
Fix $c\in \mathbb C^m \setminus \{0\}$ and take positive integers $q,N,k$ with $ q>(N-k+1)(k+1), N\geq k$. Let $V$ be a complex projective variety of dimension $k$ embedding into $\PP^n(\C)$. Let $D_j\ (1\leq j\leq q)$ be a hypersurface of degree $d_j$  in $\mathbb P^n(\C)$ such that $D_1, \dots, D_q$ are in $N$-subgeneral position on $V.$ Let $d$ be the least common multiple of
$d_1, \dots, d_q.$ Let $ f=(f_0:f_1:\dots:f_{n}): \C^m \to V$  be a algebraically non-degenerate meromorphic mapping on
$\mathcal P^1_{c}$ with $\varsigma(f)=\varsigma <1.$ Then for any $\varepsilon >0$, we have
\begin{align*}
\quad (q-N_k-\varepsilon)T_f(r)&\le \sum_{j=1}^{q}d_j^{-1}\overset{\sim}
N_{f}^{[M, c]}(r, D_j)+S(r, f),
\end{align*}
where $N_k=(N-k+1)(k+1),$ $M=\dfrac{(dk)^{k} \deg V}{k!}\Big(1+2ld^k\deg V (N-k+1)(2k+1)I(\varepsilon^{-1})\Big)^{k}$, in which $I(x):=\min\{k\in \N: k>x\}$ for $x>0.$
\end{theorem}

Recall that the hypersurfaces $D_1, \dots, D_q$ are said to be {\it in $N$-subgeneral position on V} if for every subset
$\{i_1, \dots, i_{N+1}\} \subset \{1, \dots, q\},$ we have
$$ V\cap \text{Supp}D_{i_{1}} \cap \dots \cap \text{Supp}D_{i_{N+1}}= \emptyset,$$
where $\text{Supp}(D)$ means the support of the divisor $D.$ If $N=k,$ the hypersurfaces $D_1, \dots, D_q$ are said to be {\it in general position on V}. In Theorem~\ref{th1}, $\mathcal P^{1}_{c}$ is the field of meromorphic functions $h$ of period $c$ in $\mathbb C^m$ with the hyper-order $\varsigma(h)<1$. Suppose that $V$ is generated by a homogeneous ideal
$\mathcal I(V )$ in $\mathbb{C}[x_0,\dots, x_n] $ and let
$\mathcal I_{\mathcal P^1_{c}}(V)$ be the  ideal in $\mathcal P^1_{c}[x_0,\dots, x_n] $ generated by $\mathcal I(V ),$ so that $P(f_0, \dots, f_n)=0$ for all $P\in \mathcal I_{\mathcal P^1_{c}}(V).$
The mapping $f$ is said to be {\it algebraically nondegenerate} over $\mathcal P^1_{c}$ if there is not $Q\in \mathcal P^1_{c}[x_0,\dots, x_n] \setminus \mathcal I_{\mathcal P^1_{c}}(V)$  such that $Q(f_0, \dots, f_n)=0.$ When $V=\mathbb P^n(\mathbb C),$ $f$ is algebraically nondegenerate over $\mathcal P^1_{c}$ if there is not a nonzero homogeneous polynomial $Q\in \mathcal P^1_{c}[x_0,\dots, x_n] $ such that $Q(f_0, \dots, f_n)=0.$

Theorem \ref{th1} is a difference analogue of second main theorem due to S. D. Quang \cite{Q}, which  implies the following defect relation immediately:
\begin{equation*}
 \sum_{j=1}^{q}\overset{\sim} \delta_f^{[M, c]}(r, D_j)\le (N-k+1)(k+1)+\varepsilon,
\end{equation*}
where $\overset{\sim} \delta_f^{[M, c]}(r, D_j)=1-\limsup_{r\to \infty}\dfrac{\overset{\sim}
N_{f}^{[M, c]}(r, D_j)}{d_jT_f(r)}.$ When $N=k,$ it  is a difference analogue of the defect relation due to M. Ru \cite{Ru3}.

Next we show that by using Veronese imbedding, the second main theorem can be modified such that the step number $M$ of difference does not depend on $\varepsilon.$
Let $D$ be a hypersurface of degree $d$ in $\mathbb P^n(\mathbb C)$, which is defined by a homogeneous polynomial $Q$ of degree $d$
$$Q({\bf x}) = \sum\limits_{j=0}^{n_d}a_j{\bf x}^{I_j},$$
where ${\bf x}=(x_0,\dots,x_n)\in \mathbb{C}^{n+1}$, $I_j=(i_{j0},\dots ,i_{jn}) \in\mathbb{Z}_+^{n+1}$ is a $(n+1)$-fold index over the set $\mathbb{Z}_+$ of nonnegative integers with $|I_j|=i_{j0}+\dots +i_{jn} =d$ for $j=0, \dots, n_d$, $n_d= \binom{n+d}{n}-1$, $\ab = (a_0,\dots,a_{n_d})\in \mathbb{C}^{n_d+1}$, ${\bf x}^{I_j}=x_0^{i_{j0}}\dots x_n^{i_{jn}}$, and where $I_0, \dots,I_{n_d}$ are {\it lexicographic ordering}. Let $w=(w_0:\dots:w_{n_d})$ be a homogeneous coordinate in $\PP^{n_d}(\mathbb C)$ and let  $\varrho_{d}: \mathbb P^n(\mathbb C) \to \mathbb P^{n_d}(\mathbb C)$ be the {\it Veronese embedding of degree
$d$} defined by
\begin{equation*}
    \varrho_{d}(x)= (w_0(x):\dots:w_{n_{\mathcal D}}(x)) = ({\bf x}^{I_0}:\dots:{\bf x}^{I_{n_d}}),
\end{equation*}
where $x=(x_0:\dots:x_{n})$ is a homogeneous coordinate in $\PP^{n}(\mathbb C)$. Then a linear form
\begin{equation*}
    L(w) =  a_{0}w_0 + \dots + a_{{n_d}}w_{n_d}
\end{equation*}
and a hyperplane $H=\{w:L(w) =0\}$ in $\mathbb P^{n_d}(\mathbb C)$ are associated with $D$ (or $Q$). By a way, Casorati determinant of the mapping $f: \mathbb C^m \to \mathbb P^n(\mathbb C)$ respect to Veronese embedding $\varrho_{d}$ and ${c}\in \C^m\setminus \{0\}$ is defined by
$$ C_{\varrho_{d}}(f(z)) =
\left|\begin{array}{cccc}
{\tilde{f}}^{I_0}(z)&{\tilde{f}}^{I_1}(z)&\cdots & {\tilde{f}}^{I_{n_d}}(z)\\
{\tilde{f}}^{I_0}(z+{c})&{\tilde{f}}^{I_1}(z+{ c})&\cdots & {\tilde{f}}^{I_{n_d}}(z+{c})\\
\vdots&\vdots&\ddots&\vdots\\
{\tilde{f}}^{I_0}(z+n_d{c})&{\tilde{f}}^{I_1}(z+n_d{ c})&\cdots & {\tilde{f}}^{I_{n_d}}(z+n_d{ c})\\
\end{array}
\right|.$$
When $d=1,$ we call
$$C(f(z)):= C_{\varrho_{1}}(f(z))=\left|\begin{array}{cccc}
f_{0}(z)&f_{1}(z)&\cdots & f_{n}(z)\\
f_{0}(z+{c})&f_{1}(z+{ c})&\cdots & f_{n}(z+{c})\\
\vdots&\vdots&\ddots&\vdots\\
f_{0}(z+n{c})&f_{1}(z+n{ c})&\cdots & f_{n}(z+n{ c})\\
\end{array}
\right|$$  by the Casorati determinant of $f.$

Let $\mathcal D = \{D_1,\dots,D_q\}$ be a collection of arbitrary hypersurfaces and let $Q_j$ be the homogeneous polynomial of degree $d_j$ in $\mathbb C[x_0, \dots, x_n]$ defining $D_j$ for $ j =1,\dots,q$. Let $d=d_{\mathcal D}$ be the least common multiple of $d_1,...,d_q$ and set
$$n_{\mathcal D}=n_d= \binom{n+d}{n}-1. $$
Set $Q^*_j = Q_j^{d/d_j}$, so that a vector $\ab^*_j\in \mathbb{C}^{n_{\mathcal D}+1}$ and a hyperplane $H_j$ in $\mathbb P^{n_{\mathcal D}}(\mathbb C)$ are associated with $Q^*_j$ for $j=1,\dots, q.$ If $q> n_{\mathcal D}$, $N\geq n_{\mathcal D}$, the collection $\mathcal D$ is said to be in {\it $N$-subgeneral position for Veronese embedding} in $\mathbb P^n(\mathbb C)$ if $\{H_1,\dots,H_q\} $ is in $N$-subgeneral position in $\mathbb P^{n_{\mathcal D}}(\mathbb C)$, that is,  for any distinct indices
 $i_0,\dots,i_{N} $ lied in $ \{1,\dots,q\}$, the vectors $\ab^*_{i_0},\dots,\ab^*_{i_{N}}$ have rank $n_{\mathcal D}.$

\begin{theorem}\label{th2a}
Let $ f=(f_0:f_1:\dots:f_{n}): \mathbb C^m \to \mathbb P^n(\mathbb C)$  be algebraically non-degenerate meromorphic mapping. Take an integer $N$ with $N\geq n_{\mathcal D}$ and let $D_j$ $(1\leq j\leq q)$ be a hypersurface of degree $d_j$ in $\mathbb P^n(\mathbb C)$ such that $D_1, \dots, D_q$ are in $N$-subgeneral position  for Veronese embedding.  Let $d$ be the least common multiple of
$d_1, \dots, d_q$. If $ q>2N-n_{\mathcal D}+1$,  then for any $1< r <+\infty,$
\begin{align*}
  \quad (q-(2N-n_{\mathcal D}+1))T_f(r)&\le \sum_{j=1}^{q}d_j^{-1} N_{f}^{n_{\mathcal D}}(r, D_j)+S(r, f).
\end{align*}
\end{theorem}

\begin{theorem}\label{th3a}
Take an integer $N$ with $N\geq n_{\mathcal D}$ and let $D_j\ (1\leq j\leq q)$ be a hypersurface of degree $d_j$ in $\mathbb P^n(\C)$ such that $D_1, \dots, D_q$ are in $N$-subgeneral
position for Veronese embedding.  Let $d$ be the least common multiple of
$d_1, \dots, d_q$ and set $d'=\min_{j=1}^{q}\{d_j\}.$  Let $ f, g : \mathbb C^m \to \mathbb P^n(\mathbb C)$  be algebraically non-degenerate meromorphic mappings satisfying the following conditions:
\begin{description}
  \item[(a)] $\dim f^{-1}( D_{i})\cap f^{-1}(D_j)\le m-2$ and $\dim g^{-1}(D_{i})\cap g^{-1}(D_j)\le m-2$ for all $1\le i<j\le q;$
  \item[(b)] there exist positive integers $m_j$ with $m_1\ge \dots \ge m_q$ such that
\begin{equation*}
    \overline E_{m_j)}(D_j, f)= \overline E_{m_j)}(D_j, g),\ j=1, \dots, q,
\end{equation*}
and $f(z)=g(z)$ for $z\in\cup_{j=1}^{q}\overline E_{m_j)}(D_j, f)$.
\end{description}
Then $f\equiv g$ if
$$ q>2N-n_{\mathcal D}+1+\sum_{i=1}^{q}\dfrac{n_{\mathcal D}}{m_i+1}+\dfrac{2m_1n_{\mathcal D}}{d'(m_1+1)}.$$
\end{theorem}

If we take $N=n_{\mathcal D}$ and let $m_j\to\infty$ for each $j$ in Theorem \ref{th3a}, we get a uniqueness theorem under $q\geq 3n_{\mathcal D}+2.$ The number of hypersurfaces in Theorem \ref{th2a} is smaller than that appeared in the result of Dulock-Ru \cite{DR}.

\begin{theorem}\label{th2}
Let $D_j\ (1\leq j\leq q)$ be a hypersurface of degree $d_j$ in $\mathbb P^n(\C)$ such that $D_1, \dots, D_q$ are in $N$-subgeneral position for Veronese embedding.  Let $d$ be the least common multiple of
$d_1, \dots, d_q$. Let $ f=(f_0:f_1:\dots:f_{n}): \mathbb C^m \to \mathbb P^n(\mathbb C)$  be algebraically non-degenerate meromorphic mapping on
$\mathcal P^1_{c}$ with $\varsigma(f) <1.$ If $q>2N-n_{\mathcal D}+1,$ $N\geq n_{\mathcal D}$, we have
\begin{align*}
 \quad (q-(2N-n_{\mathcal D}+1))T_f(r)&\le \sum_{j=1}^{q}d_j^{-1}N_f(r, D_j)-\dfrac{N}{n_{\mathcal D}}N_{C_{\varrho_{d}}(f)}(r)+S(r, f)\\
&\le \sum_{j=1}^{q}d_j^{-1}\overset{\sim}
N_{f}^{[n_{\mathcal D}, c]}(r, D_j)+S(r, f).
\end{align*}
\end{theorem}

When $d=1,$ then $n_D=n,$ Theorem~\ref{th2} is just the difference analogue of the second main theorem of meromorphic mappings due to Cao-Korhonen \cite{TC}.

\begin{cor}\label{cor2}
Take $q=n_{\mathcal D}+2$ and let $D_j\ (1\leq j\leq q)$ be a hypersurface of degree $d_j$ in $\mathbb P^n(\C)$ such that $D_1, \dots, D_q$ are in $N$-subgeneral position for Veronese embedding.
Let $d$ be the least common multiple of
$d_1, \dots, d_q.$ Then a nonconstant meromorphic mapping $ f=(f_0:f_1:\dots:f_{n}) : \mathbb C^m \to \mathbb P^n(\mathbb C)$ is algebraically degenerate on $\mathcal P^1_{c}$ if  $\varsigma(f) <1$ and if
\begin{align*}
  \sum_{j=1}^{q}d_j^{-1}\overset{\sim} N_{f}^{[n_{\mathcal D}, c]}(r, D_j)=S(r, f).
\end{align*}
\end{cor}

If the preimages of hypersurfaces $D_j\ (j=1, \dots, n_{\mathcal D}+2)$ under $f$ are forward invariant with respect to the translation $\tau(z)=z+c,$ that is, $\tau(f^{-1}(D))\subset f^{-1}(D)$ counting multiplicity,
then  Corollary \ref{cor2} means that $f$ is algebraically degenerate on $\mathcal P^1_{c}.$

\begin{theorem}\label{th3}
Take $c\in \mathbb C^m\setminus \{0\}$ and let $f: \mathbb C^m\to \mathbb P^n(\mathbb C)$ be a meromorphic mapping with $\varsigma(f) <1.$  Let $D_j\ (j=1, \dots, q)$ be a hypersurface of degree $d_j$ in $P^n(\mathbb C)$ such that $D_1,...,D_q$ are in $N$-subgeneral position for Veronese embedding. Let $d$ be the least common multiple of
$d_1, \dots, d_q.$ If the preimages of hypersurfaces $D_j$ under $f$ are forward invariant with respect to the translation $\tau(z)=z+c,$ then the image of $f$ is contained in one of hypersurfaces $D_j$ or the image of $\varrho_{d}(f)$ is contained in a projective linear subspace of dimension $\le \left[\dfrac{N}{q-N}\right]$ over $\mathcal P_c^{1}$. In particular, we have $f= f_c:=f\circ\tau$ if $q\ge 2N+1$.
\end{theorem}

When $d=1,$ then $n_D=n.$ If we choose $N=n, q=n+p$ in Theorem \ref{th3},  it is a difference analogue of Picard's theorem for holomorphic curves due to Halburd-Korhonen-Tohge \cite{RRK}.

Take $p=(p_1, \dots, p_m)\in \mathbb C^m \setminus \{0, 1\}$ and write $p\cdot z=(p_1z_1, \dots, p_mz_m)$ for $z=(z_1, \dots, z_m)$. Let $\mathcal P^{0}_{p}$ be the field of meromorphic functions $h$  of zero order in $\mathbb C^m$ such that $h(p\cdot z)=h(z).$ Similarly, the meromorphic mapping $f: \mathbb C^m \to V\subset \mathbb P^n(\mathbb C)$ is said to be {\it algebraically nondegenerate} over $\mathcal P^0_{p}$ if there is not $Q\in \mathcal P_{p}^{0}[x_0, \dots, x_n]\setminus \mathcal I_{\mathcal P^0_{p}}(V)$ such that $Q(f_0, \dots, f_n)=0$, where $\mathcal I_{\mathcal P^0_{p}}(V)$ is the  ideal in $\mathcal P^0_{p}[x_0,\dots, x_n] $ generated by $\mathcal I(V ).$ When $V=\mathbb P^n(\mathbb C),$ $f$ is algebraically nondegenerate over $\mathcal P^0_{p}$ if there is not a polynomial homogeneous $Q\in \mathcal P^0_{p}[x_0,\dots, x_n] \setminus \{0\}$ such that $Q(f_0, \dots, f_n)=0.$ By a way, the $p$-Casorati determinant of a meromorphic mapping $f:\C^m \to \mathbb P^n(\C)$ is defined by
\[C_p(f(z)) =
\left|\begin{array}{cccc}
f_{0}(z)&f_{1}(z)&\cdots & f_{n}(z)\\
f_{0}(p\cdot z)&f_{1}(p\cdot z)&\cdots & f_{n}(p\cdot z)\\
\vdots&\vdots&\ddots&\vdots\\
f_{0}(p^n\cdot z)&f_{1}(p^n\cdot z)&\cdots & f_{n}(p^n\cdot z)\\
\end{array}
\right|,
\]
in which $p^j=(p_1^j, \dots, p_m^j)$. By a result in \cite{TC1}, we know that $C_p(f)\equiv 0$ if and only if $f_0, \dots, f_n$ are linear dependent over the filed $\mathcal P^0_{p}.$
By using the $p$-difference analogue of logarithmic derivative lemma (cf.\cite{LDP} and \cite{TC1}), we can obtain the following results  for $p$-difference operator without proofs.

\begin{theorem}\label{th1p}
Take positive integers $N,k$ with $N\geq k$ and let $V\subset \mathbb P^{n}(\mathbb C)$ be a complex projective variety of dimension $k\ge 1.$ Let $D_j\ (1\leq j\leq q)$
 be a hypersurface  of degree $d_j$ in $\mathbb P^n(\mathbb C)$ such that $D_1, \dots, D_q$ are in $N$-subgeneral position on $V.$ Let $d$ be the least common multiple of
$d_1, \dots, d_q.$ Let $ f=(f_0:f_1:\dots:f_{n}):\mathbb C^m \to V$  be algebraically non-degenerate meromorphic mapping of zero order on
$\mathcal P^0_{p}$. Assume that $q>N_k.$ Then for any $\varepsilon >0,$ we have
\begin{align*}
 (q-N_k-\varepsilon)T_f(r)&\le \sum_{j=1}^{q}d_j^{-1}\hat N_{f}^{[M, p]}(r, D_j)+S_1(r, f),
\end{align*}
where
$M=\dfrac{(dk)^{k}\deg V}{k!}\Big(1+2ld^k\deg V (N-k+1)(2k+1)I(\varepsilon^{-1})\Big)^{k}.$
\end{theorem}

\begin{theorem}\label{th2p}
Let $ f=(f_0:f_1:\dots:f_{n}):\mathbb C^m \to \mathbb P^n(\mathbb C)$  be algebraically non-degenerate meromorphic mapping of zero order on
$\mathcal P^0_{p}$. Let $D_j\ (1\leq j\leq q)$ be a hypersurface of degree $d_j$  in $\mathbb P^n(\C)$ such that $D_1, \dots, D_q$ are in $N$-subgeneral position for Veronese embedding. Let $d$ be the least common multiple of
$d_1, \dots, d_q$. If $q>2N-n_{\mathcal D}+1,$ $N\geq n_{\mathcal D}$, then we have
\begin{align*}
 \quad (q-(2N-n_{\mathcal D}+1))T_f(r)&\le \sum_{j=1}^{q}d_j^{-1}\hat N_{f}^{[n_{\mathcal D}, p]}(r, D_j)+S_1(r, f).
\end{align*}
\end{theorem}

\begin{theorem}\label{th3p}
Let $f=(f_0:f_1:\dots:f_{n}):\mathbb C^m\to \mathbb P^n(\mathbb C)$ be a meromorphic mapping of zero order. Let $D_j\ ( j=1, \dots, q)$ be a hypersurface of degree $d_j$ in $ \mathbb P^n(\mathbb C)$ such that $D_1,\dots,D_q$ are in $N$-subgeneral position for Veronese embedding. Let $d$ be the least common multiple of
$d_1, \dots, d_q.$ If  the preimages of hypersurfaces $D_j$ under $f$  are forward invariant with respect to the rescaling $\tau(z)=p\cdot z,$ then the image of $f$ is contained in one of hypersurfaces $D_j$ or the image of
 $\varrho_{d}(f)$ is contained in a projective linear subspace of dimension $\le \left[\dfrac{N}{q-N}\right]$ over $\mathcal P^0_{p}$. In particular, we have $f= f\circ\tau$ if $q\ge 2N+1$.
\end{theorem}

\section{Some Lemmas}
\def\theequation{2.\arabic{equation}}
\setcounter{equation}{0}

In order to prove theorems above, we need the following lemmas.

\begin{lemma}\label{lm21}\cite{RRK, TC}
Let $f$ be a non-constant meromorphic function in $\mathbb C^m,$ and  let ${c}\in \C^m\setminus \{0\}.$ If $
\varsigma(f)<1,$ then the function $f_c(z)=f(z+c)$ satisfies
$$ m\left(r, \dfrac{f_c}{f}\right) =S(r, f).$$
\end{lemma}

From Lemma \ref{lm21}, if we repalce $f$ by $\dfrac{1}{f},$ and using First Main Theorem,
we have $\varsigma\left(\dfrac{1}{f}\right)=\varsigma(f)<1,$ then
\begin{equation}\label{rm1}
  m\left(r, \dfrac{f}{f_c}\right) =S(r, f).
\end{equation}

\begin{lemma}\label{lm23}\cite{RRK, TC}
Take $c\in \mathbb C^m\setminus\{0\}$. If a meromorphic mapping $g=[g_0:\dots :g_n]:\mathbb{C}^m\longrightarrow \mathbb{P}^n(\mathbb C)$ satisfies $\varsigma(g)<1$,
then the Casorati determinant of $g$ satisfies $C(g)\equiv 0$ if and only if the entire functions $g_0, \dots, g_n$ are linearly dependent
 over the field $\mathcal P^{1}_{c}.$
\end{lemma}

\begin{lemma}\label{lm23b}\cite{RRK, TC} Take $c\in \mathbb C^m\setminus \{0\}$ and let $g=(g_0: \dots: g_n):\mathbb{C}^m\longrightarrow \mathbb{P}^n(\mathbb C)$ be a meromorphic mapping with $\varsigma(g)<1$ such that  all zeros of $g_0, \dots, g_n$ are forward invariant with respect to the translation $\tau(z)=z+c.$ Let
$ S_1\cup \dots \cup S_l $
 be a partition of $\{0, \dots, n\}$ formed in such a way that $i$ and $j$ are the same class $S_k$ if and only if $g_i/g_j\in \mathcal P_{c}^{1}.$ If
$ g_0+\dots+g_n=0, $
then for each $k\in\{1,\dots,l\}$, we have
$$ \sum_{i\in S_k}g_i=0. $$
\end{lemma}

\begin{lemma}\label{lm23a}\cite{RRK}
Let $T : [0, +\infty) \to [0, +\infty)$ be a non-decreasing continuous function
such that the hyper-order of $T$ is strictly less than one, that is,
$$\varsigma=\limsup_{r \to \infty}\dfrac{\log \log T(r)}{\log r}<1.$$
If $\delta \in (0, 1 -\varsigma ),$ $s \in (0, \infty)$, then
$$ T(r+s)=T(r)+o\Big(\dfrac{T(r)}{r^{\delta}}\Big), $$
where $r$ runs to infinity outside of a set of finite logarithmic measure.
\end{lemma}

Take ${\bf c}=(c_0, \dots, c_n)\in \mathbb{R}^{n+1}$. Let $V\subset \mathbb P^n(\mathbb C)$ be a projective variety of dimension $k$ and let
$$ F_V({\bf x}_{0}, \dots, {\bf x}_{k})=F_V(x_{00}, \dots, x_{0n}; \dots; x_{k0}, \dots, x_{kn}) $$
be the Chow form associated to $V$, that is, $F_V$ is irreducible in $\mathbb C[x_{00}, \dots, x_{kn}],$ $F_V$ is a homogeneous polynomial of degree $\Delta=\deg V$ in each block ${\bf x}_i=(x_{i0}, \dots, x_{in})$ $ (0\leq i\leq k)$ and $F_V({\bf x}_0, \dots, {\bf x}_k)=0$ if and only if $X\cap H_{{\bf x}_0}\cap \dots\cap H_{{\bf x}_k}\ne \emptyset,$ where $H_{{\bf x}_i}$ is a hyperplane given by
$$H_{{\bf x}_i}=\{(u_0:\cdots :u_n)\ | \  x_{i0}u_0+\dots+u_{in}u_n=0\}. $$
For an auxiliary variable $t$, we obtain a decomposition
\begin{align*}
F_V(t^{c_0}x_{00},& \dots, t^{c_n}x_{0n};\dots;t^{c_0}x_{k0}, \dots, t^{c_n}x_{kn})\\
&=t^{e_0}G_0({\bf x}_0, \dots, {\bf x}_k)+\dots+t^{e_r}G_r({\bf x}_0, \dots, {\bf x}_k),
\end{align*}
where $G_0, \dots, G_r\in \mathbb C[x_{00},\dots, x_{0n};\dots; x_{k0}, \dots, x_{kn}]$ and $e_0>e_1>\dots>e_r.$ The {\it Chow weight} of $V$ with respect to ${\bf c}$ is defined by
$$ e_V({\bf c}):=e_0. $$

Suppose that $V$ is generated by a prime ideal
$\mathcal I(V )$ in $\mathbb{C}[x_0,\dots, x_n] $ and let $\C[x_0, \dots, x_n]_u$ be the vector space of homogeneous polynomials in $\C[x_0, \dots, x_n]$ of degree $u$ (including 0). Put $\mathcal I_u(V ):=\C[x_0, \dots, x_n]_u\cap \mathcal I(V )$ and denote the Hilbert function $H_V$ of $V$ by
$$ H_V(u):=\dim V_u,\; V_u= \C[x_0, \dots, x_n]_u / \mathcal I_u(V).$$
We define the $u$-th Hilbert weight
$S_V(u, {\bf c})$ of $V$ with respect to ${\bf c}$ by
$$ S_V(u, {\bf c}):=\max \left(\sum_{j=1}^{H_V(u)}I_j\cdot {\bf c}\right),$$
where the maximum is taken over all sets of monomials ${\bf x}^{I_1}, \dots, {\bf x}^{I_{H_V(u)}}$ whose residue classes module $\mathcal I(V )$ form a basis of $\mathbb C[x_0, \dots, x_n]_u/\mathcal I_u(V)$, in which $I_j=(i_{j0},\dots ,i_{jn}) \in \mathbb{Z}_+^{n+1}$ is a $(n+1)$-fold index.

\begin{lemma}\label{lm24}\cite{Ru3}
Take ${\bf c}=(c_0, \dots, c_n) \in \mathbb R^{n+1}$ with $c_j\geq 0$ for $j=0,...,n$.
Let $V\subset \mathbb P^{n}(\mathbb C)$ be an algebraic variety of dimension $k$.   If $u>\Delta=\deg V$, then
$$ \dfrac{1}{uH_V(u)}S_V(u, {\bf c})\ge \dfrac{1}{(k+1)\Delta}e_V({\bf c})-\dfrac{(2k+1)\Delta}{u}\max_{0\le i\le n}c_i.$$
\end{lemma}

\begin{lemma}\label{lm25}\cite{EF, Ru3}
Take ${\bf c}=(c_0, \dots, c_n) \in \mathbb R^{n+1}$ with $c_j\ge 0$ for $j=0,\dots,n$. Let $V$ be a subvariety of $\mathbb P^{n}(\C)$ of dimension $k$ and let $\{i_0, \dots, i_{k}\}$ be a subset of
$\{0, \dots , n\}$ such that $\{ x_{i_0}=\dots=x_{i_{k}}=0\}\cap V= \emptyset$.  Then
$$ e_V({\bf c}) \ge (c_{i_0}+\dots + c_{i_{k}})\deg V.$$
\end{lemma}

Apply to Lemma \ref{lm21} to Lemma \ref{lm23a} and  using the idea in \cite{TC1}, we get
a result as follows.

\begin{lemma}\label{lm31} Let $ f=(f_0:f_{1}:\dots:f_{n}) : \mathbb C^m \ax \mathbb P^{n}(\mathbb C)$  be a linearly non-degenerate meromorphic mapping on $\mathcal P^{1}_{c}$ with $\varsigma(f)<1.$ Let $H_{{\bf a}_1},\dots , H_{{\bf a}_q}$ be hyperplanes in $\mathbb P^{n}(\mathbb C)$ in general position which are defined by linear forms $L_{{\bf a}_1},\dots, L_{{\bf a}_q}$ respectively. Then we have the inequality
\begin{align}\label{ba-ieq}
 \quad \int_{S(r)}\max_{K}\sum_{l \in K}\log \dfrac {\|\tilde{f}\|}{|L_{{\bf a}_l}(\tilde{f})|}\sigma \leqslant (n+1)
T_{f}(r) - N_{C(f)}(r) +S(r, f),
\end{align}
in which the maximum is taken over all subsets $K$ of $\{1,\dots,q\}$ such that $\{{\bf a}_{l}\}_ {l\in K}$ are linearly independent.
\end{lemma}

\begin{proof}
Let $K \subset \{1,\dots,q\}$ be a set such that $\{{\bf a}_{l}\}_ {l\in K}$ are linearly independent. Note that $ q\geqslant n+1$ and the cardinal number  $\#K$ of the set $K$ satisfies $\#K\leq n+1$. First of all, we consider the case $\# K=n+1.$ Let $\mathcal T$ be the set of all injective maps
$\mu:\{ 0,1,\dots,n\} \ax \{1,\dots,q\}$ and set $h_{l}=L_{{\bf a}_{l}}(\tilde{f})$. The quantity of left side at the inequality (\ref{ba-ieq}), say $A(r)$ simply, satisfies
\begin{align*}
A(r) &=\int_{S(r)}\max_{\mu \in \mathcal  T}\sum_{l=0}^n\log \dfrac {\|\tilde{f}\|}{|h_{{\mu(l)}}|}
\sigma +O(1)\\
&=\int_{S(r)}\max_{\mu \in \mathcal  T}\log \left\{ \dfrac {\|\tilde{f}\|^{n+1}}{\prod\limits_{l=0}^n|h_{{\mu(l)}}|}\right\} \sigma +O(1)\\
&\le \int_{S(r)}\max_{\mu \in \mathcal  T} \log \bigg\{ \dfrac {\|\tilde{f}\|^{n+1}}
{|C(h_{{\mu(0)}},\dots,h_{{\mu(n)}})|}\bigg\}\sigma\\
&+ \int_{S(r)}\max_{\mu \in \mathcal  T}\log \bigg\{ \dfrac {|C(h_{{\mu(0)}},\dots, h_{{\mu(n)}})|}
{\prod\limits_{l=0}^n|h_{{\mu(l)}}|}\bigg\}\sigma+O(1),
\end{align*}
which further implies
\begin{align*}
A(r) &\le\int_{S(r)}\log \bigg\{\max_{\mu \in \mathcal  T}\dfrac {\|\tilde{f}\|^{n+1}}
{|C(h_{{\mu(0)}},\dots, h_{{\mu(n)}})|}\bigg\}\sigma +B(r)+O(1)\\
&\le\int_{S(r)}\log \bigg\{\sum_{\mu \in \mathcal  T}\dfrac {\|\tilde{f}\|^{n+1}}
{|C(h_{{\mu(0)}},\dots, h_{{\mu(n)}})|}\bigg\}\sigma +B(r)+O(1),
\end{align*}
where
$$ B(r)=\sum_{\mu \in T} \int_{S(r)}\log  \bigg\{\dfrac {|C(h_{{\mu(0)}},\dots, h_{{\mu(n)}})|}
{\prod\limits_{l=0}^n|h_{{\mu(l)}}|}\bigg\}\sigma.$$
By using Lemma \ref{lm23}, we easily find the following relation of Casorati determinants
$$|C(h_{{\mu(0)}},\dots,h_{{\mu(n)}})| =C_{\mu}|C(f)|\not\equiv 0,$$
where $C_{\mu} \ne 0$ is a constant. Thus we obtain
\begin{equation}\label{m1}
A(r)\leq \int_{S(r)}\log \dfrac{\|\tilde{f}\|^{n+1}} {|C(f)|}\sigma +B(r)+O(1).
\end{equation}

Next we estimate $B(r)$. Let  $\mathbb Z^{+}$ be the set of positive integers and set
\begin{equation*}
 h_{\mu(j),tc}(z)= h_{{\mu(j)}}(z+tc),\ g_{\mu(j)}=\dfrac{h_{\mu(j)}}{h_{\mu(0)}}
\end{equation*}
for $t\in \mathbb Z^{+}$, $j\in\{0, \dots, n\}$.
By a simple calculation, we find
$$ C(h_{\mu(0)}, \dots, h_{\mu(n)})=\sum_{\theta}\text{sgn}(\theta)h_{\theta(\mu(0))}h_{\theta(\mu(1)),c}\dots h_{\theta(\mu(n)),nc},$$
where the sum runs over all permutations  $$\{\theta: \{\mu(0), \dots , \mu(n)\} \to  \{\mu(0), \dots, \mu(n)\}\}$$ of $n+1$ objects, which further implies
\begin{align*}
&\dfrac{C(h_{\mu(0)}, \dots, h_{\mu(n)})}{h_{\mu(0)}h_{\mu(1),c} \dots h_{\mu(n),nc}}\\
&=  \dfrac{1}{ g_{\mu(1),c}\dots g_{\mu(n),nc}}\sum_{\theta}\text{sgn}(\theta)g_{\theta(\mu(0))} g_{\theta(\mu(1)),c}\dots  g_{\theta(\mu(n)),nc},
\end{align*}
and hence
\begin{align}\label{ctma1}
\dfrac{|C(h_{\mu(0)}, \dots, h_{\mu(n)})|}{|h_{\mu(0)} h_{\mu(1),c} \dots  h_{\mu(n),nc}|}
\le \dfrac{\sum_{i_1+i_2+\dots+i_n\le \dfrac{n(n+1)}{2}}
\prod_{l=1}^n\dfrac{|g_{\mu(l),i_lc}|}{|g_{\mu(l)}|}}
{\dfrac{| g_{\mu(1),c}|}{|g_{\mu(1)}|} \dots \dfrac{| g_{\mu(n),nc}|}{|g_{\mu(n)}|}}.
\end{align}
Applying Jensen's formula to  meromorphic function $h_{\mu}=\dfrac{h_{\mu(1),c}\dots  h_{\mu(n),nc}}{h_{\mu(1)}\dots h_{\mu(n)}}$ on $\mathbb{C}^m$, we have
\begin{align}\label{t1}
\int\limits_{S(r)}\log|h_{\mu}|\sigma=N_{h_{\mu(1),c}\dots  h_{\mu(n),nc}}(r)-N_{h_{\mu(1)}\dots h_{\mu(n)}}(r)+O(1).
\end{align}
Note that
\begin{align}\label{t2}
 N_{h_{\mu(1),c}\dots  h_{\mu(n),nc}}(r)=\sum_{j=1}^{n}N_{h_{\mu(j), jc}}(r)
\end{align}
and
\begin{align}\label{t3}
 N_{h_{\mu(1)}\dots h_{\mu(n)}}(r)=\sum_{j=1}^{n}N_{h_{\mu(j)}}(r).
\end{align}
Combining (\ref{t1}) to (\ref{t3}), we obtain
\begin{align}\label{ctma2}
\int_{S(r)}\log|h_{\mu}|\sigma= \sum_{j=1}^{n}(N_{h_{\mu(j),jc}}(r)-N_{h_{\mu(j)}}(r))+O(1).
\end{align}
Further, note that $z_0$ is a zero of $h_{\mu(j),jc}$ if and only if $z_0+jc$ is a zero of $h_{\mu(j)}.$ Thus we have
$$N_{h_{\mu(j),jc}}(r)\le N_{h_{\mu(j)}}(r+j||c||)$$
for each $j=1, \dots, n.$  By using First Main Theorem
$$N_{h_{\mu(j)}}(r)\le T_{f}(r)+O(1)$$
and noting that
$$ \varsigma_j=\limsup_{r\to \infty}\dfrac{\log\log N_{h_{\mu(j)}}(r)}{\log r}\le \varsigma(f)<1,$$
then Lemma \ref{lm23a} deduces
\begin{align*}
 N_{h_{\mu(j)}}(r+j||c||) \le N_{h_{\mu(j)}}(r)+S(r,f)
\end{align*}
for each $j=1, \dots, n$, which immediately yields
\begin{align}\label{ctm}
N_{h_{\mu(j),jc}}(r)\le N_{h_{\mu(j)}}(r)+S(r,f),
\end{align}
and combine (\ref{ctma2}) and (\ref{ctm}), we get
\begin{align*}
\int_{S(r)}\log|h_{\mu}|\sigma\leq S(r,f).
\end{align*}
Combining above equality with (\ref{ctma1}) and using Lemma \ref{lm21} , we get
\begin{align}\label{ctb2}
B(r)&= \sum_{\mu \in \mathcal T}\int_{S(r)}\log \Big\{\dfrac{|C(h_{\mu(0)}, h_{\mu(1)}, \dots, h_{\mu(n)})|}
{|h_{\mu(0)} h_{\mu(1),c}\dots  h_{\mu(n),nc}|}\Big\}\sigma\notag
+\sum_{\mu\in \mathcal T}\int_{S(r)}\log|h_{\mu}|\sigma\\
&\le \sum_{\mu\in \mathcal T}\sum_{i_1+\dots+i_n\le \dfrac{n(n+1)}{2}}\sum_{l=1}^{n}\int_{S(r)}\log^{+}
\dfrac{| g_{\mu(l),i_lc}|}{|g_{\mu(l)}|}\sigma\notag\\
&+\sum_{\mu\in \mathcal T}\sum_{l=1}^{n}\log^{+}\dfrac{|g_{\mu(l)}|}
{|g_{\mu(l),lc}|}\sigma+S(r,f)\leq S(r,f).
\end{align}
By using Jensen's formula and combining (\ref{m1}) with (\ref{ctb2}), we obtain
\begin{align}\label{ctb3}
\int_{S(r)}\max_{K}\sum_{l \in K}\log \dfrac {\|\tilde{f}\|}{|h_{l}(\tilde{f})|}\sigma
 \leqslant (n+1)
T_{f}(r) - N_{C(f)}(r) +S(r, f).
\end{align}

Finally, we consider the case $ \# K<n+1.$ We can take $n+1- \#K$ elements from the set $\{1, \dots, q\}$ united with the set $K$ to form a new set $K'$ such that $\{{\bf a}_{l}\}_{l\in K'}$ are linearly independent.
Note that
\begin{align*}
\int_{S(r)}\max_{K}\sum_{l \in K}\log \dfrac {\|\tilde{f}\|}{|h_{l}(\tilde{f})|}\sigma
 \leqslant \int_{S(r)}\max_{K'}\sum_{l' \in K'}\log \dfrac {\|\tilde{f}\|}{|h_{{l'}}(\tilde{f})|}\sigma+O(1).
\end{align*}
We also obtain the estimate (\ref{ctb3}). The proof of Lemma \ref{lm31} is completed.
\end{proof}

\begin{lemma}\label{lmt1}\cite{TC}
Take $c\in \mathbb C^m \setminus \{0\}$ and let $f : \mathbb C^m \to \mathbb P^n(\mathbb C)$ be a linearly nondegenerate meromorphic mapping over $\mathcal P_c^{1}$ with hyperorder $\varsigma (f)< 1$. Let $H_j (1\le j\le q)$ be $q(>2N - n + 1)$ hyperplanes in $N$-subgeneral position in $\mathbb P^n(\mathbb C).$ Then we have
\begin{align*}
(q - 2N + n -1)T_f(r) &\le \sum_{j=1}^{q}N_f(r, H_j)-\dfrac{N}{n}N_{C(f)}(r)+S(r, f)\\
 &\le \sum_{j=1}^{q}\overset{\sim} N_f^{[n, c]}(r, H_j)+S(r, f).
\end{align*}
\end{lemma}

\begin{lemma}\label{lmt2}\cite{NK, N}
Let $f:\mathbb C^m\to \mathbb P^n(\mathbb C)$ be a linearly nondegenerate meromorphic mapping and let $\{H_j\}_{j=1}^{q}$ be hyperplanes in $N$-subgeneral position in $\mathbb P^n(\mathbb C)$ with $N\ge n$ and $q>2N-n+1$. Then we have
$$ (q-2N+n-1)T_f(r)\le \sum_{j=1}^{q}N_f^{n}(r, H_j)+S(r, f).$$
\end{lemma}

\begin{lemma}\label{lmt3}\cite{Q}
Take positive integers $N$ and $k$ with $N\geq k$. Let $V$ be a smooth projective subvariety of dimension $k$ in $\mathbb P^n(\mathbb C)$ and let $D_j\ (1\leq j\leq N+1)$ be a hypersurface in $\mathbb P^n(\mathbb C)$ defined by a homogenous polynomial $Q_j$ of same degree $d$  such that
$$ \big(\cap_{j=1}^{N+1}{\rm Supp}D_j\big)\cap V=\emptyset. $$
Then there exist constants $b_{tj}\in \mathbb C$ such that hypersurfaces $D^*_t\ (1\leq t\leq k+1)$ defined by $P_t$ satisfy $\big(\cap_{t=1}^{k+1}{\rm Supp}D^*_t\big)\cap V=\emptyset,$ in which
\begin{equation*}
    P_t=\left\{
          \begin{array}{ll}
            Q_1, & \hbox{ $t=1$;} \\
            \sum_{j=2}^{N-k+t}b_{tj}Q_j, & \hbox{$t\geq 2$.}
          \end{array}
        \right.
\end{equation*}
\end{lemma}

\section{Proof of Theorem \ref{th1}}

Suppose that the ideal $\mathcal I(V)$ is generated by homogeneous polynomials $\mathcal P_1$, ..., $\mathcal P_{\alpha}$ and let $\mathcal H_j\ (1\leq j\leq \alpha)$ be a  hypersurface defined by $\mathcal P_j$. Let $Q_i\ (1\le i\le q)$ be a homogeneous polynomial of degree $d_i$ in $\C[x_0, \dots, x_n]$ defining the hypersurface $D_i.$ Since we can replace $Q_i$ by $Q_i^{d/d_i},$  where $d$ is the least common multiple of $d_1, \dots, d_q$, we may assume that $Q_1, \dots, Q_q$ have the same degree $d$. By the assumption, $D_1,\dots,D_q$ are in $N$-subgeneral position on $V$.  Then for every subset
$\mathfrak{i}=\{i_1, \dots, i_{N+1}\} \subset \{1, \dots, q\},$ we have
$$ V\cap \text{Supp}D_{i_{1}} \cap \dots \cap \text{Supp}D_{i_{N+1}}= \emptyset,$$
which means
$$ \text{Supp}{\mathcal H_1}\cap\dots\cap\text{Supp}{\mathcal H_{\alpha}}\cap \text{Supp}D_{i_{1}} \cap \dots \cap \text{Supp}D_{i_{N+1}}= \emptyset.$$
Thus by Hilbert's Nullstellensatz \cite{VW}, for each integer $\beta \in \{0, \dots, n\},$ there is an integer $m_{\beta} > \{d, \max_{t=1}^{\alpha}\{\deg \mathcal P_t\}\}$ such that
$$ x_{\beta}^{m_{\beta}}=\sum_{j=1}^{N+1}a_{{\beta}j}(x_0, \dots, x_n) Q_{i_j}(x_0, \dots, x_n)
+\sum_{t=1}^{\alpha}b_{\beta t}(x_0, \dots, x_n)\mathcal P_t(x_0, \dots, x_n),$$
where $a_{{\beta}j}$ (resp. $b_{\beta t}$) is a homogeneous polynomial of degree $m_{\beta}-d$ (resp. $m_{\beta}-\deg P_t$) in $\mathbb{C}[x_0,\dots,x_n]$. Since $f(\mathbb{C}^m)\subseteq V$, we have
$$ \sum_{t=1}^{\alpha}b_{\beta t}(f_0(z), \dots, f_n(z))\mathcal P_t(f_0(z), \dots, f_n(z))=0,\ z\in \mathbb{C}^m.$$
Thus we obtain an estimate
$$
|f_{\beta}(z )|^{m_\beta}\le \kappa_{\mathfrak i} ||\tilde{f}(z)||^{m_{\beta}-d}\max_{1\leq j\leq N+1}\{|Q_{i_j}(\tilde{f}(z))|\},$$
where $\kappa_{\mathfrak i}$ is a positive constant depends only on the coefficients of $a_{{\beta}j}$ $(0\le \beta \le n, 1\le j\le N+1)$. Set $A=\max_{\mathfrak i} \kappa_{\mathfrak i}$.  We have
\begin{align}\label{2a}
 ||\tilde{f}(z)||^{d}\le A\max_{1\le j\le N+1}|Q_{i_j}(\tilde{f}(z))|,\ z\in \mathbb{C}^m.
\end{align}

Let $\mathbb I$ be the set of all permutations of the set $\{1, \dots, q\}$ with the cardinal number $n_0:=q!$, which can be listed as $\mathbb I=\{ \mathbb I_1, \dots, \mathbb I_{n_0}\},$ where $\mathbb I_i=(\mathbb I_i(1), \dots, \mathbb I_i(q))\in (\mathbb Z^+)^{q}$ and $\mathbb I_1<\dots<\mathbb I_{n_0}$ according to the lexicographic order.
For each $i\in\{1,\dots,n_0\}$, Lemm \ref{lmt3} implies that there exist constants $b_{itj}\in \mathbb C$ such that hypersurfaces $D^*_{i,t}\ (1\leq t\leq k+1)$ defined by $P_{i,t}$ satisfy $\big(\cap_{t=1}^{k+1}{\rm Supp}D^*_{i,t}\big)\cap V=\emptyset,$ in which
\begin{equation}\label{ctb1}
    P_{i,t}=\left\{
          \begin{array}{ll}
            Q_{\mathbb{I}_i(1)}, & \hbox{ $t=1$;} \\
            \sum_{j=2}^{N-k+t}b_{itj}Q_{\mathbb{I}_i(j)}, & \hbox{$t\geq 2$.}
          \end{array}
        \right.
\end{equation}
Then there exists a positive constant $B\ge 1$ such that
\begin{align}\label{tm1}
 |P_{i, t}({\bf x})|\le B\max_{1\le j\le N-k+t}|Q_{\mathbb{I}_i(j)}({\bf x})|
\end{align}
for $1\le t\le k+1$, ${\bf x}=(x_0, \dots, x_n)\in \mathbb C^{n+1}.$ Moreover, there exists a positive constant $E$ such that
$$ |P_{i, t}({\bf x})|\le E||{\bf x}||^{d}.$$

Further, one defines a mapping $\widetilde{\psi}:\mathbb{C}^{n+1}\longrightarrow \mathbb{C}^l$ by
\begin{equation*}
 \widetilde{\psi}({\bf x})=   (P_{1, 1}({\bf x}),\dots,P_{1, k+1}({\bf x}),  \dots,P_{n_0, 1}({\bf x}), \dots ,P_{n_0,k+1}({\bf x})),
\end{equation*}
where $l=n_0(k+1)$, and let ${\mathfrak e}=({\mathfrak e}_0:{\mathfrak e}_1:\cdots:{\mathfrak e}_n):V\longrightarrow \mathbb{P}^n(\mathbb{C})$ be the embedding mapping. Then a finite morphism $\psi:V\longrightarrow\mathbb P^{l-1}(\mathbb C)$ is well defined by
\begin{equation*}
 \psi(x)=   (P_{1, 1}({\bf x}):\dots:P_{1, k+1}({\bf x}):  \dots:P_{n_0, 1}({\bf x}): \dots : P_{n_0,k+1}({\bf x})),
\end{equation*}
where ${\bf x}=({\mathfrak e}_0(x),...,{\mathfrak e}_n(x))$, such that $Y=\psi(V)$ is a complex projective subvarieties of $\PP^{l-1}(\C)$ with $\dim Y=k,$ $\deg Y:=\Delta\le d^k \deg V$ (cf. \cite{MS}).

Moreover, taking a positive integer $u$ and fixed a basis $\{\phi_0, \dots, \phi_{n_u}\}$ of  the vector space $Y_u:={\mathbb C[y_1, \dots, y_l]_u}/ {{\mathcal I}_u(Y)}$, where $n_u+1=H_Y(u)=\dim Y_u$, $\mathcal I_u(Y):=\C[y_1, \dots, y_l]_u\cap \mathcal I(Y)$ in which $\mathcal I(Y)$ is the prime ideal which defines algebraic variety $Y,$ according to S. D. Quang \cite{Q}, a meromorphic mapping
$$ F=(\phi_0(\psi(f)):\dots:\phi_{n_u}(\psi(f))): \mathbb C^m \to \mathbb P^{n_u}(\mathbb C)$$
is well defined with a reduced representation $\tilde{F}=\left(\dfrac{\phi_0(\tilde{\psi}(\tilde{f}))}{\varphi},\dots,\dfrac{\phi_{n_u}(\tilde{\psi}(\tilde{f}))}{\varphi}\right),$ where $\varphi$ is a common factor of $\phi_0(\psi(f)),\dots,\phi_{n_u}(\psi(f))$ which is a holomorphic function on $\mathbb C^m$. Moreover, $F$ is linearly non-degenerate on $\mathcal P^1_{c}$. Note that
\begin{align*}
T_{F}(r)&=\int\limits_{S(r)}\log \sqrt{\dfrac{\sum_{j=0}^{n_u}|\phi_j(\tilde\psi(\tilde f(z)))|^2}{\varphi^2(z)}}\sigma(z)+O(1)\notag\\
&=\int\limits_{S(r)}\log \sqrt{\sum_{j=0}^{n_u}|\phi_j(\tilde\psi(\tilde f(z)))|^2}\sigma(z)-\int\limits_{S(r)}\log |\varphi(z)|\sigma(z)+O(1).
\end{align*}
Then we obtain an estimate
\begin{align*}
T_{F}(r)&\le du\int\limits_{S(r)}\log ||\tilde f(z)||\sigma(z)-\int\limits_{S(r)}\log |\varphi(z)|\sigma(z)+O(1)\notag\\
&=du T_f(r)-N_{\varphi}(r)+O(1)\le du T_f(r)+O(1),
\end{align*}
which implies
\begin{align}\label{ma3}
S(r, F)\le S(r, f).
\end{align}

Now, we claim
\begin{align}\label{Claim 1.1}
\log \prod_{i=1}^{q} \dfrac{||\tilde{f}(z)||^{d}}{|Q_i(\tilde{f}(z))|}\le (N-k+1)\max_{i}\log \dfrac{||\tilde{f}(z)||^{(k+1)d}}{\prod_{t=1}^{k+1}|P_{i, t}(\tilde{f}(z))|}+O(1)
\end{align}
for $z\in \mathbb{C}^m\setminus \mathfrak S$, where $\mathfrak S=\mathfrak S_1\cup \mathfrak S_2$ is defined by
\begin{equation*}
    \mathfrak S_1=\bigcup_{i=1}^{q}(Q_i(\tilde f))^{-1}(0) ,\ \mathfrak S_2=\bigcup_{1\le i\le n_0,1\le t\le k+1}(P_{i, t}(\tilde f))^{-1}(0),
\end{equation*}
so that (\ref{Claim 1.1}) yields immediately
\begin{align}\label{3aa}
\sum_{i=1}^{q}m_f(r,D_i)\le (N-k+1)\int\limits_{S(r)}\max_{i}\log \dfrac{||\tilde{f}(z)||^{(k+1)d}}{\prod_{t=1}^{k+1}|P_{i, t}(\tilde{f}(z))|}\sigma+O(1).
\end{align}
In fact, for an fixed $z\in \mathbb{C}^m\setminus \mathfrak S$, there exists some $i\in\{1,...,n_0\}$ such that
\begin{align}\label{1a}
|Q_{\mathbb I_i(1)}(\tilde{f}(z))| \le  |Q_{\mathbb I_i(2)}(\tilde{f}(z))| \le  \dots \le |Q_{\mathbb I_i(q)}(\tilde{f}(z))|.
\end{align}
By (\ref{1a}) and (\ref{2a}), we have $ \dfrac{||\tilde{f}(z)||^{d}}{|Q_{\mathbb I_i(j)}(\tilde{f}(z))|}\le A $
 for all $j\geq N+1$, which means
 \begin{equation*}
 \prod_{i=1}^{q}\dfrac{||\tilde{f}(z)||^{d}}{|Q_i(\tilde{f}(z))|}   \le A^{q-N-1}\prod_{j=1}^{N+1}\dfrac{||\tilde{f}(z)||^{d}}{|Q_{\mathbb I_i(j)}(\tilde{f}(z))|}.
 \end{equation*}
It follow from  (\ref{1a}) and (\ref{tm1}) that
$$ \prod_{j=1}^{N-k+1}|Q_{\mathbb I_i(j)}(\tilde{f}(z))|\ge |Q_{\mathbb I_i(1)}(\tilde{f}(z))|^{N-k+1}=|P_{i,1}(\tilde{f}(z))|^{N-k+1}$$
and
$$ |Q_{\mathbb I_i(N-k+t)}(\tilde{f}(z))|\ge \dfrac{|P_{i, t}(\tilde{f}(z))|}{B} ,\ t=2,\dots,k+1,$$
so that
\begin{equation*}
\prod_{i=1}^{q}\dfrac{||\tilde{f}(z)||^{d}}{|Q_i(\tilde{f}(z))|} \leq A^{q-N-1}B^k\frac{||\tilde{f}(z)||^{(N+1)d}}{|P_{i, 1}(\tilde{f}(z))|^{N-k+1}\prod_{t=2}^{k+1}|P_{i, t}(\tilde{f}(z))|}.
\end{equation*}
Since
\begin{equation*}
\left(\prod_{t=2}^{k+1}|P_{i, t}(\tilde{f}(z))|\right)^{N-k}\leq E^{(N-k)k}||\tilde{f}(z)||^{(N-k)kd},
\end{equation*}
we deduce
\begin{equation}\label{3a}
\prod_{i=1}^{q}\dfrac{||\tilde{f}(z)||^{d}}{|Q_i(\tilde{f}(z))|} \leq A^{q-N-1}B^kE^{(N-k)k}\dfrac{||\tilde{f}(z)||^{(N+1)d+(N-k)kd}}{\prod_{t=1}^{k+1}|P_{i, t}(\tilde{f}(z))|^{N-k+1}}.
\end{equation}
Note that $ (N+1)d+(N-k)kd=N_kd$, where $N_k=(N-k+1)(k+1)$. Then (\ref{3a}) yields the claim (\ref{Claim 1.1}) immediately.

Next, by using the auxiliary mapping $F:\mathbb{C}^m\longrightarrow \mathbb{P}^{n_u}(\mathbb{C})$, we claim
\begin{align}\label{Claim 1.2}
d(q-N_k)T_f(r)&\le \sum_{i=1}^{q}N_f(r, D_i)-\dfrac{N_k}{uH_Y(u)}N_{C(F)}(r)-\dfrac{N_k}{u}N_{\varphi}(r)\notag\\
 &+\dfrac{N_k(2k+1)\Delta}{u}\sum_{1\le i\le n_0, 1\le t\le k+1}m_f(r, D_{i,t}^{*})+S(r, f),
\end{align}
where $N_{C(F)}(r)$  is the valence function for the zeros of Casorati determinant $C(F)$ of $F$.
Fixed a point $z\in \mathbb{C}^m\setminus \mathfrak S$  and taking an index $i\in \{1, \dots, n_0\}$, we define
$${\bf c}_z=(c_{1,1,z},\dots, c_{1,k+1, z}, c_{2,1,z}, \dots, c_{2,k+1, z}, c_{n_0,1,z}\dots, c_{n_0, k+1,z})\in \mathbb R_{+}^{l},$$ where $c_{i,t,z}=\log \dfrac{||\tilde{f}(z)||^{d}||P_{i,t}||}{|P_{i,t}(\tilde{f}(z))|}$ for $i\in\{1, \dots, n_0\}, t\in\{1, \dots, k+1\}$, and $\mathbb R_{+}$ is the set of nonnegative real numbers.
Then according to the definition of the $u$-th Hilbert weight $S_Y(u, {\bf c}_z)$ of $Y$ with respect to ${\bf c}_z$, there exists a subset  $A_z\subset \{0, \dots, l_u\}$ with $
  |A_z|=n_u+1=H_Y(u)$ in which $ l_u=\binom{l+u-1}{u}-1$ such that $\{{\bf x}^{{I}_{j,z}}: j\in A_z\}$ is a basis of the vector space $Y_u$ and
\begin{align*}
S_Y(u, {\bf c}_z)= \sum_{j\in A_z}I_{j,z}\cdot {\bf c}_z,
\end{align*}
where $ {\bf x}\in \mathbb{C}^l$, and
\begin{equation*}
    {I}_{j,z} = (a_{j,1,1,z},\dots, a_{j,1,k+1,z},\dots, a_{j,n_0,1,z},\dots, a_{j,n_0,k+1,z})\in \Z_{+}^{l},
\end{equation*}
that is,
\begin{align*}
S_Y(u, {\bf c}_z)&=\sum_{j\in A_z}\sum_{1\le i\le n_0,1\le t\le k+1}a_{j,i,t,z}\log  \dfrac{||\tilde f(z)||^{d}||P_{i,t}||}{|P_{i,t}(\tilde f(z))|}\notag\\
&=-\log \prod_{j,i,t}|P_{i,t}(\tilde f(z))|^{a_{j,i,t,z}}+\log ||\tilde f(z)||^{duH_Y(u)}+O(uH_Y(u)),
\end{align*}
or equivalently,
\begin{equation*}
\frac{1}{uH_Y(u)}\log \prod_{j,i,t}|P_{i,t}(\tilde f(z))|^{a_{j,i,t,z}}=-\frac{S_Y(u, {\bf c}_z)}{uH_Y(u)}+\log ||\tilde f(z)||^d+O(1).
\end{equation*}
Since $\{\phi_0, \dots, \phi_{n_u}\}$ also is a  basis of $Y_u$, there exist $n_u+1$ linear forms $$\mathcal L_z=\{L_{j, z}, j\in A_z\}$$ which are linearly independent such that
\begin{align*}
{\bf x}^{{I}_{j,z}}=L_{j, z}(\phi_0({\bf x}), \dots, \phi_{n_u}({\bf x})).
\end{align*}
Note that
\begin{align*}
\widetilde{\psi}(\tilde f)^{{I}_{j,z}}=L_{j,z}(\phi_0(\widetilde{\psi}(\tilde f)),\dots, \phi_{n_u}(\widetilde{\psi}(\tilde f)))&=\varphi L_{j,z}(\tilde F).
\end{align*}
Then we obtain
$$\prod_{i, t}P_{i,t}(\tilde{f}(z))^{a_{j,i,t,z}}=\widetilde{\psi}(\tilde f(z))^{{I}_{j,z}}=\varphi (z)L_{j,z}(\tilde F(z)).$$
Thus we have
\begin{align*}
\frac{1}{uH_Y(u)}\log \prod_{j\in A_z}|\varphi (z)L_{j,z}(\tilde F(z))|
&=\frac{1}{uH_Y(u)}\log \prod_{j,i,t} |P_{i,t}(\tilde{f}(z))|^{a_{j,i,t,z}}\\
&=-\frac{S_Y(u, {\bf c}_z)}{uH_Y(u)}+\log ||\tilde{f}(z)||^d+O(1),
\end{align*}
and hence
\begin{align}\label{mb2}
\frac{1}{uH_Y(u)} \log \prod_{j\in A_z}&\dfrac{||\tilde{F}(z)||.||L_{j,z}||}{|L_{j,z}(\tilde F(z))|}=\frac{S_Y(u, {\bf c}_z)}{uH_Y(u)}-\log ||\tilde{f}(z)||^d\notag\\
&+\frac{1}{u}\log ||\tilde{F}(z)||+\frac{1}{u}\log |\varphi (z)|+O(1).
\end{align}

Further, by using Lemma \ref{lm24} and Lemma \ref{lm25}, we obtain easily
\begin{eqnarray}\label{ct44a1}
\frac{S_Y(u,{\bf c}_z)}{uH_Y(u)} &\geq& \dfrac{1}{k+1}\max_i\log \dfrac{||\tilde{f}(z)||^{(k+1)d}}{\prod_{t=1}^{k+1}|P_{i,t}(\tilde{f}(z))|} \nonumber\\
 && - \frac{(2k+1)\Delta}{u} \sum_{i,t}\log \dfrac{||\tilde{f}(z)||^{d}||P_{i,t}||}{|P_{i,t}(\tilde{f}(z))|}
+O(1).
\end{eqnarray}
Note that when $z$ is changed, the subset $A_z$ of $\{0,...,l_u\}$ runs only over a finite set, so that  $\mathcal L=\cup_{z}\mathcal L_z$ also is a finite set with $|\mathcal L|\le \binom{l_u+1}{n_u+1}$.
Combine (\ref{mb2}) and (\ref{ct44a1}), we get
\begin{align}\label{tmc}
&\frac{1}{k+1}\int\limits_{S(r)}\max_{i}\log \dfrac{||\tilde{f}(z)||^{(k+1)d}}{\prod_{t=1}^{k+1}|P_{i,t}(\tilde{f}(z))|}
\le dT_f(r)-\frac{1}{u}T_F(r)-\dfrac{1}{u}N_{\varphi}(r)\notag\\
&+\dfrac{1}{uH_Y(u)} \int\limits_{S(r)}\max_{J\subset \mathcal L} \log \prod_{L\in J}\dfrac{||\tilde{F}(z)||.||L||}{|L(\tilde F(z))|}\sigma\notag\\
&+\dfrac{(2k+1)\Delta}{u}\sum_{i, t}m_f(r,D_{i,t}^{*})
+O(1),
\end{align}
where $J$ is a subset of $\mathcal L$ with $|J|=H_Y(u)$ such that $\{L: L\in J\}$ is linearly independent. Applying Lemma \ref{lm31} to $F$, we see
\begin{align}\label{ct49a}
\int_{S(r)}\max_{J\subset \mathcal L} \log \prod_{L\in J}\dfrac{||\tilde{F}(z)||.||L||}{|L(\tilde{F}(z))|}\sigma\le (n_u+1)T_F(r)-N_{C(F)}(r)+S(r, F).
\end{align}
Thus by (\ref{3aa}) and (\ref{tmc}), we have
\begin{align}\label{ct46a}
\sum_{i=1}^{q}m_f(r,D_i)&\le dN_kT_f(r)-\frac{N_k}{uH_Y(u)}N_{C(F)}(r)-\dfrac{N_k}{u}N_{\varphi}(r)\notag\\
&+\dfrac{N_k(2k+1)\Delta}{u}\sum_{i, t}m_f(r,D_{i,t}^{*})
+S(r,F).
\end{align}
Hence the claim (\ref{Claim 1.2})  follows from (\ref{ma3}) and  First Main Theorem immediately.

Thirdly, we claim
\begin{align}\label{Claim:1.3}
\sum_{i=1}^{q}&N_f(r, D_i)-\dfrac{N_k}{uH_Y(u)}N_{C(F)}(r)-\dfrac{N_k}{u}N_{\varphi}(r)
\le \sum_{i=1}^{q} \overset{\sim}N_{f}^{[n_u, c]}(r, D_i)\notag\\
&+\dfrac{N_k(2k+1)\Delta}{u} \sum_{1\le i\le n_0, 1\le t\le k+1} N_f(r, D_{i,t}^{*}),
\end{align}
which will follows from a functional inequality
\begin{align}\label{tcm}
\sum_{i=1}^{q}&\nu_{Q_i(\tilde f)}-\dfrac{N_k}{uH_Y(u)}\nu_{C(F)}-\dfrac{N_k}{u}\nu_{\varphi}
\le \sum_{i=1}^{q} {\overset{\sim}\nu}^{[n_u, c]}_{Q_i(\tilde f)}\notag\\
&+\dfrac{N_k(2k+1)\Delta}{u} \sum_{1\le i\le n_0, 1\le t\le k+1} \nu_{P_{i,t}(\tilde f)}.
\end{align}
In fact, (\ref{tcm}) is trivial over $\mathbb C^m\setminus \mathfrak S_1$. Obviously, we only need to prove (\ref{tcm}) for $z\in \mathfrak S_1\backslash I_f.$
Write
\begin{equation*}
\mathfrak S_1 \backslash I_f=\bigcup_{i=1}^{n_0}    \mathfrak T(i),
\end{equation*}
where $\mathfrak T(i)$ is the set of all points $z\in \mathfrak S_1\backslash I_f$ satisfying (\ref{1a}).
Take $z\in \mathfrak S_1\backslash I_f$. W.l.o.g., we may assume $z\in \mathfrak T(1)$, where $\mathbb I_1=(1, \dots, q).$ Since $D_1, \dots, D_q$ are in $N$-subgeneral position in $V,$ then there are at most $N$ hypersurfaces in $\{D_1, \dots, D_q\}$ intersecting with $V.$ Thus (\ref{1a}) with $i=1$ implies that there exists $\mathfrak p \in \{1,\dots,N\}$ such that
$\nu_{Q_j(\tilde f)}(z)>0$ for all $j=1,\dots,\mathfrak p,$ and $\nu_{Q_j(\tilde f)}(z)=0$ for all $j=\mathfrak p+1,\dots,q.$ Hence by renumbering the set $\{1, \dots, q\}$ if necessary, we may assume that
\begin{align}\label{ctb1a}
\nu_{Q_1(\tilde f)}(z)\ge \dots\ge \nu_{Q_{\mathfrak p}(\tilde f)}(z)>0=\nu_{Q_{\mathfrak p+1}(\tilde f)}(z)=\dots=\nu_{Q_{q}(\tilde f)}(z).
\end{align}
Thus we further claim
\begin{align}\label{h}
(N-k+1)\sum_{t=1}^{k+1}\nu_{P_{1,t}(\tilde{f})}^{[n_u,c]}(z)\ge\sum_{i=1}^{q}\nu_{Q_i(\tilde{f})}^{[n_u,c]}(z).
\end{align}
In fact, for the case $\mathfrak p=1$, (\ref{h}) follows easily from
\begin{align}\label{hb}
(N-k+1)\sum_{t=1}^{k+1}\nu_{P_{1,t}(\tilde{f})}^{[n_u,c]}(z)\ge\nu_{Q_1(\tilde{f})}^{[n_u,c]}(z)=\sum_{i=1}^{q}\nu_{Q_i(\tilde{f})}^{[n_u,c]}(z).
\end{align}
When $\mathfrak p\ge 2$, we next prove (\ref{h}) by distinguishing two cases.

{\bf Case 1.} $k=1.$

For this case, there is only one hypersurface in $D_{1,1}^{*},\dots, D_{1,k+1}^{*}$ intersecting with $V$, which is just $D_{1,1}^{*}$ since  $\nu_{P_{1,1}(\tilde f)}(z)=\nu_{Q_1(\tilde f)}(z)>0$.
Thus there are at most $N-k$ functions in $\{Q_{2}(\tilde f),\dots,Q_{N-k+2}(\tilde f)\}$ vanishing at $z,$ that is, $N-k+2>\mathfrak p.$ Indeed, if all these functions vanish at $z,$ then $\nu_{P_{1,2}(\tilde f)}(z)>0$, that is, $D_{1,2}^{*}$ intersects with $V$. This is a contradiction.
Hence  (\ref{h}) follows from
\begin{equation}\label{hb0}
(N-k+1)\sum_{t=1}^{k+1}\nu_{P_{1,t}(\tilde{f})}^{[n_u,c]}(z)\geq \mathfrak p\nu_{Q_1(\tilde{f})}^{[n_u,c]}(z)\geq    \sum_{i=1}^{\mathfrak p}\nu_{Q_i(\tilde{f})}^{[n_u,c]}(z) =\sum_{i=1}^{q}\nu_{Q_i(\tilde{f})}^{[n_u,c]}(z).
\end{equation}

{\bf Case 2.} $k\ge 2.$

Moreover, if $N-k+2>\mathfrak p,$ then the estimate (\ref{hb0}) still holds, so that it is sufficient to consider the case $N-k+2\leq \mathfrak p$. Thus there exists $\mathfrak l\in\{2,k+1\}$ such that $N-k+\mathfrak l=\mathfrak p.$ By the definition of $P_{1,t}$ at (\ref{ctb1}), we see easily
\begin{equation}\label{hh1}
   \nu_{P_{1,t}(\tilde{f})}(z)\ge \min_{2\leq j\leq N-k+t} \nu_{Q_j(\tilde{f})}(z)=\nu_{Q_{N-k+t}(\tilde f)}(z)>0
\end{equation}
for all $t=2,\dots,\mathfrak l$. Since there exists $\mathfrak h\in \{1,\dots, k\}$ such that $\mathfrak h$ hypersurfaces in $D_{1,1}^{*},\dots, D_{1,k+1}^{*}$ intersect with $V$,  (\ref{hh1}) means $2\le \mathfrak l\le \mathfrak h\le k.$ Note that any $n_u$-successive and $c$-separated point $z$ of $Q_{N-k+t}(\tilde{f})\ (t=2,\dots,\mathfrak l)$ is a $n_u$-successive and $c$-separated point of $P_{1,t}(\tilde f)\ ( 2\le t\le \mathfrak l)$, that is, $\nu_{Q_{N-k+t}(\tilde f)}^{[n_u,c]}(z)=\nu_{Q_{N-k+t}(\tilde f)}(z)$  and $ \nu_{P_{1,t}(\tilde{f})}^{[n_u,c]}(z)= \nu_{P_{1,t}(\tilde{f})}(z)$ for $t=2,\dots,\mathfrak l.$  We see
\begin{align*}
\sum_{t=1}^{k+1}&\nu_{P_{1,t}(\tilde{f})}^{[n_u,c]}(z)=\sum_{t=1}^{\mathfrak h}\nu_{P_{1,t}(\tilde{f})}^{[n_u,c]}(z)\notag
\ge \nu_{Q_1(\tilde{f})}^{[n_u,c]}(z)+\sum_{t=2}^{\mathfrak l}\nu_{P_{1,t}(\tilde{f})}^{[n_u,c]}(z),
\end{align*}
which implies
\begin{align}\label{hd}
(N-k+1)\sum_{t=1}^{k+1}&\nu_{P_{1,t}(\tilde{f})}^{[n_u,c]}(z)\ge (N-k+1) \nu_{Q_1(\tilde{f})}^{[n_u,c]}(z)+\sum_{t=2}^{\mathfrak l}\nu_{Q_{N-k+t}(\tilde{f})}^{[n_u,c]}(z)\notag\\
&\ge \sum_{i=1}^{N-k+\mathfrak l}\nu_{Q_i(\tilde{f})}^{[n_u,c]}(z)= \sum_{i=1}^{\mathfrak p}\nu_{Q_i(\tilde{f})}^{[n_u,c]}(z)=\sum_{i=1}^{q}\nu_{Q_i(\tilde{f})}^{[n_u,c]}(z).
\end{align}
Hence the claim (\ref{h}) is proved.

According to the definition of $u$-th Hilbert weight $S_Y(u, {\bf c})$ with respect to
$$ {\bf c}=(c_{1,1},\dots, c_{1,k+1},\dots,c_{n_0,1},\dots,c_{n_0,k+1})\in\mathbb R_{+}^{l}$$
with $c_{i,t}=\nu_{P_{i,t}(\tilde{f})}^{[n_u,c]}(z)$,  there exist $H_Y(u)$ multi-indies
\begin{equation*}
  {I}_j=(a_{j,1,1},\dots,a_{j,1,k+1},\dots,a_{j,n_0,1},\dots,a_{j,n_0,k+1})\in\mathbb Z_{+}^{l}
\end{equation*}
such that ${\bf x}^{{I}_1},\dots, {\bf x}^{{I}_{H_Y(u)}}$ is a basis of the vector space $Y_u$ for $ {\bf x}\in \mathbb{C}^l$ and
\begin{align*}
S_Y(u, {\bf c})= \sum_{j=1}^{H_Y(u)}I_{j}\cdot {\bf c}.
\end{align*}
Moreover, there exist $H_Y(u)$ linear forms $L_1,\dots, L_{H_Y(u)}$ which are linearly independent such that
$${\bf x}^{{I}_{j}}=L_{j}(\phi_0({\bf x}), \dots, \phi_{n_u}({\bf x})),$$
which implies
\begin{align*}
\prod_{1\le i\le n_0, 1\le t\le k+1}P_{i,t}(\tilde{f}(z))^{a_{j,i,t}}=\widetilde{\psi}(\tilde f(z))^{{I}_{j}}=\varphi (z)L_{j}(\tilde F(z)),
\end{align*}
so that
\begin{equation}\label{b2}
\nu_{L_j(\tilde{F})}^{[n_u,c]}(z)+\nu_{\varphi}^{[n_u,c]}(z)=\sum_{i,t}a_{j,i,t}\nu_{P_{i,t}(\tilde{f})}^{[n_u,c]}(z)={I}_j\cdot {\bf c}.
\end{equation}
Note that there exists a constant $C\ne 0$ such that
$$C(L_1(\tilde{F}),\dots, L_{H_Y(u)}(\tilde{F}))=CC(F)$$
and hence
\begin{align}\label{b1}
\nu_{C(F)}^{[n_u, c]}(z)=\nu_{C(L_1(\tilde{F}),\dots,L_{H_Y(u)}(\tilde{F}))}^{[n_u,c]}(z)\ge \sum_{j=1}^{H_Y(u)}\nu_{L_j(\tilde{F})}^{[n_u,c]}(z).
\end{align}
Combining (\ref{b1}) and (\ref{b2}), we obtain
$$ \nu_{C(F)}^{[n_u,c]}(z)\ge S_Y(u,{\bf c})-H_Y(u)\nu_{\varphi}^{[n_u,c]}(z).$$
Since $D_{1,1}^{*}\dots, D_{1,k+1}^{*}$ are in general position in $V,$ Lemma \ref{lm24} and Lemma \ref{lm25} yield
\begin{align*}
S_Y(u, &{\bf c})\ge \dfrac{uH_Y(u)}{(k+1)\Delta}e_Y({\bf c})-H_Y(u)(2k+1)\Delta \max_{i, t}c_{i,t}\notag\\
&\ge \dfrac{uH_Y(u)}{(k+1)}\sum_{t=1}^{k+1}\nu_{P_{1,t}(\tilde{f})}^{[n_u,c]}(z)-H_Y(u)(2k+1)\Delta \max_{i, t}\nu_{P_{i,t}(\tilde{f})}^{[n_u,c]}(z).
\end{align*}
Thus we have
\begin{align}\label{e}
 &\dfrac{k+1}{uH_Y(u)}\nu_{C(F)}^{[n_u,c]}(z)\ge \sum_{t=1}^{k+1}\nu_{P_{1,t}(\tilde{f})}^{[n_u,c]}(z)-\dfrac{k+1}{u}\nu_{\varphi}^{[n_u,c]}(z)\notag\\
&-\dfrac{(k+1)(2k+1)\Delta}{u}\max_{i, t}\nu_{P_{i,t}(\tilde{f})}^{[n_u,c]}(z).
\end{align}
It follows from (\ref{h}) and (\ref{e}) that
\begin{align*}
 \dfrac{N_k}{uH_Y(u)}\nu_{C(F)}^{[n_u,c]}(z)
&\ge \sum_{i=1}^{q}\nu_{Q_i(\tilde{f})}^{[n_u,c]}(z)-\dfrac{N_k}{u}\nu_{\varphi}^{[n_u,c]}(z)\notag\\
&-\dfrac{N_k(2k+1)\Delta}{u}\max_{1\le i\le n_0, 1\le t\le k+1}\nu_{P_{i,t}(\tilde{f})}^{[n_u,c]}(z),\\
\end{align*}
which means
\begin{align}\label{d}
 &\sum_{i=1}^{q}\nu_{Q_i(\tilde{f})}^{[n_u,c]}(z)-\dfrac{N_k}{uH_Y(u)}\nu_{C(F)}(z)\le \dfrac{N_k}{u}\nu_{\varphi}(z)\notag\\
&+\dfrac{N_k(2k+1)\Delta}{u}\sum_{1\le i\le n_0, 1\le t\le k+1}\nu_{P_{i,t}(\tilde{f})}(z)
\end{align}
since $\nu_{C(F)}(z)\geq \nu_{C(F)}^{[n_u,c]}(z)$, $\nu_{\varphi}^{[n_u,c]}(z)\le \nu_{\varphi}(z)$ and $\nu_{P_{i,t}(\tilde{f})}^{[n_u,c]}(z)\le\nu_{P_{i,t}(\tilde{f})}(z)$. Thus the claim (\ref{tcm}) follows from (\ref{d}) and the fact
\begin{align*}
\nu_{Q_i(\tilde{f})}^{[n_u,c]}(z)=\nu_{Q_i(\tilde{f})}(z)-{\overset{\sim}\nu}_{Q_i(\tilde{f})}^{[n_u,c]}(z).
\end{align*}

Finally, the claims  (\ref{Claim 1.2}), (\ref{Claim:1.3}) and First Main Theorem for $D_{i,t}^{*}$ yield
\begin{align}\label{ct51a}
d(q-N_k)T_f(r)\le \sum_{i=1}^{q}\overset{\sim} N_{f}^{[n_u, c]}(r, D_i)+\dfrac{N_k(2k+1)\Delta}{u}ldT_f(r)+S(r, f).
\end{align}
Now, for any $\varepsilon>0,$ we choose  $u=N_k(2k+1)\Delta l I(\varepsilon^{-1})$, so that
\begin{align}\label{ct53a}
\dfrac{N_k(2k+1)\Delta l}{u}<\varepsilon
\end{align}
holds. Thus it follows from (\ref{ct51a}) and (\ref{ct53a}) that
\begin{align}\label{c1}
 d(q-N_k-\varepsilon)T_f(r)&\le \sum_{i=1}^{q}\overset{\sim} N_{f}^{[n_u, c]}(r, D_i)+S(r, f).
\end{align}

Generally, replacing $D_i$ by $D_i^{d/d_i}$ with $d={\rm lcm} \{d_1, \dots, d_q\}$, we see that
$$ \overset{\sim} N_{f}^{[n_u, c]}(r, D_i^{d/d_i})\le \dfrac{d}{d_i}\overset{\sim} N_{f}^{[n_u, c]}(r, D_i).$$
Hence (\ref{c1}) implies
\begin{align*}
 (q-N_k-\varepsilon)T_f(r)&\le \sum_{i=1}^{q}d_i^{-1} \overset{\sim} N_{f}^{[n_u, c]}(r, D_i)+S(r, f).
\end{align*}
Note that $\deg Y= \Delta \le d^k \deg V,$ $\dim Y =k$ and  
\begin{equation*}
n_u\le \Delta \binom{k+u}{k}   < \Delta \Big(1+\dfrac{u}{k}\Big)^{k}\dfrac{k^k}{k!}
\end{equation*}
(cf. \cite{MS}).
For the choice of $u,$ we have
$$ n_u\le \dfrac{(kd)^{k}\deg V}{k!}\Big(1+2ld^k\deg V (N-k+1)(2k+1)I(\varepsilon^{-1})\Big)^{k},$$
and hence the proof of Theorem \ref{th1} is completed.

\section{Proof Theorem \ref{th2a}}

Applying Lemma \ref{lmt2} to the meromorphic mapping
$$G=\varrho_{d}(f)=(\tilde{f}^{I_0}:\dots:\tilde{f}^{I_{n_{\mathcal D}}}):\mathbb C^m \to \mathbb P^{n_{\mathcal D}}(\mathbb C)$$
and the collection of hyperplanes $\mathfrak H =\{H_1,\dots,H_q\}$  in $N$-subgeneral position in $\mathbb P^{n_{\mathcal D}}(\mathbb C),$ we have
$$ (q-(2N-n_{\mathcal D}+1))T_G(r)\le \sum_{j=1} ^{q}N_G^{n_{\mathcal D}}(r, H_j)+S(r, G),$$
so that Theorem \ref{th2a} follows from the facts $ T_G(r)=d T_f(r)$ and
$$N_G^{n_{\mathcal D}}(r, H_j)=\dfrac{d}{d_j}N_f^{n_{\mathcal D}}(r, D_j).$$

\section{Proof of Theorem \ref{th3a}}

Assume, to the contrary, that $f\not\equiv g.$ Then there exists two distinct indices $s, t\in \{0, \dots, n\}$  such that $\Phi:=f_sg_t-f_tg_s\not\equiv 0$.
The conditions ${\bf (a)}$ and ${\bf (b)}$ means that points in
$\cup_{j=1}^{q}\overline E_{m_j)}(D_j, f)\backslash I_f$ are zeros of $\Phi$, so that
\begin{equation}\label{a41}
\sum_{j=1}^{q}(N_{f,\le m_j}^{1}(r, D_j)+N_{g,\le m_j}^{1}(r, D_j))\le 2N_{\Phi}(r)
\le 2(T_f(r)+T(r,g))+O(1).
\end{equation}
Note that
\begin{equation*}
N_{f, \ge m_i+1}^{n_{\mathcal D}}(r, D_i) \leq \dfrac{n_{\mathcal D}}{m_i+1}N_{f,\ge m_i+1}(r, D_i)   ,
\end{equation*}
\begin{eqnarray*}
N_{f, \le m_i}^{n_{\mathcal D}}(r, D_i) &=& \dfrac{m_i}{m_i+1}N_{f, \le m_i}^{n_{\mathcal D}}(r, D_i) + \dfrac{1}{m_i+1}N_{f, \le m_i}^{n_{\mathcal D}}(r, D_i) \\
&\leq& \dfrac{n_{\mathcal D}m_i}{m_i+1}N_{f, \le m_i}^{1}(r, D_i) + \dfrac{n_{\mathcal D}}{m_i+1}N_{f, \le m_i}(r, D_i),
\end{eqnarray*}
that is
\begin{eqnarray*}
N_f^{n_{\mathcal D}}(r, D_i)   &=&  N_{f, \le m_i}^{n_{\mathcal D}}(r, D_i)+N_{f, \ge m_i+1}^{n_{\mathcal D}}(r, D_i) \\
  &\leq &  \dfrac{n_{\mathcal D}m_i}{m_i+1}N_{f, \le m_i}^{1}(r, D_i) + \dfrac{n_{\mathcal D}}{m_i+1}N_{f,\ge m_i+1}(r, D_i).
\end{eqnarray*}
Applying Theorem \ref{th2a} and First Main Theorem, we have
\begin{align*}
(q-(2N&-n_{\mathcal D}+1))T_f(r)\le \sum_{i=1}^{q}\dfrac{1}{d_i}N_f^{n_{\mathcal D}}(r, D_i)+S(r, f)\notag\\
&\le \sum_{i=1}^{q} \dfrac{n_{\mathcal D}m_i}{d_i(m_i+1)}N_{f, \le m_i}^{1}(r, D_i)+\sum_{i=1}^{q}\dfrac{n_{\mathcal D}}{m_i+1}T_f(r) +S(r, f),
\end{align*}
which implies
\begin{equation}\label{a6}
q'T_f(r)\leq \dfrac{m_1n_{\mathcal D}}{d'(m_1+1)}\sum_{j=1}^{q}N_{f, \le m_i}^{1}(r, D_i)+S(r, f),
\end{equation}
where
\begin{equation*}
    q'=q-(2N-n_{\mathcal D}+1)-\sum_{i=1}^{q}\dfrac{n_{\mathcal D}}{m_i+1}.
\end{equation*}
Similarly, we get
\begin{equation}\label{a7}
    q'T_g(r)\le \dfrac{m_1n_{\mathcal D}}{d'(m_1+1)}\sum_{j=1}^{q}N_{g, \le m_i}^{1}(r, D_i)+S(r, g).
\end{equation}
Combining (\ref{a41}), (\ref{a6}) and (\ref{a7}), we have
\begin{equation*}
\left(q'-\dfrac{2m_1n_{\mathcal D}}{d'(m_1+1)}\right)(T(r, f)+T(r,g))\leq  S(r, f)+S(r, g),
\end{equation*}
which means $q'-\dfrac{2m_1n_{\mathcal D}}{d'(m_1+1)}\leq 0$.
This is  contradiction.
Then $f\equiv g$ and the proof of Theorem \ref{th3a} is completed.

\section{Proof of Theorem \ref{th2}}

Applying Lemma \ref{lmt1} to the meromorphic mapping
$$G=\varrho_{d}(f)=(\tilde{f}^{I_0}:\dots:\tilde{f}^{I_{n_{\mathcal D}}}):\mathbb C^m \to \mathbb P^{n_{\mathcal D}}(\mathbb C)$$
and the collection of hyperplanes $\mathfrak H =\{H_1,\dots,H_q\}$  in $N$-subgeneral position in $\mathbb P^{n_{\mathcal D}}(\mathbb C),$ we have
\begin{align*}
 (q-(2N-n_{\mathcal D}+1))T_G(r)&\le \sum_{j=1}^{q} N_G(r, H_j)-\dfrac{N}{n_{\mathcal D}}N_{C(G)}(r)+S(r, G)\notag\\
&\le \sum_{j=1} ^{q}\overset{\sim} N_G^{[n_{\mathcal D},c]}(r, H_j)+S(r, G).
\end{align*}
Thus Theorem \ref{th2} follows from $ T_G(r)=dT_f(r)$ and the fact
$$\overset{\sim} N_G^{[n_{\mathcal D},c]}(r, H_j)=\dfrac{d}{d_j}\overset{\sim} N_f^{[n_{\mathcal D},c]}(r, D_j).$$

\section{Proof of Theorem \ref{th3}}

We consider the collection of hyperplanes $\mathfrak H =\{H_1,\dots,H_q\}$ associated with $\{Q_1^{*}, \dots, Q_q^{*}\}$ in $N$-subgeneral position in $\mathbb P^{n_{\mathcal D}}(\mathbb C)$, which are defined by
$$H_j =\{w=(w_0:\dots:w_{n_{\mathcal D}})\in \mathbb P^{n_{\mathcal D}}(\mathbb C): a_{j0}w_0 + \dots + a_{j{n_{\mathcal D}}}w_{n_{\mathcal D}}=0\}$$
and use  the meromorphic mapping
$$G=\varrho_{d}(f):\mathbb C^m \to \mathbb P^{n_{\mathcal D}}(\mathbb C)$$
with a reduced representation $\tilde{G}=(\tilde{f}^{I_0},\dots,\tilde{f}^{I_{n_{\mathcal D}}})$. Then
$$ G_j=H_j(\tilde{G})=\sum_{i=0}^{n_{\mathcal D}}a_{ji}{\tilde{f}}^{I_{i}}$$
satisfies $\tau({G_j}^{-1}(0))\subset {G_j}^{-1}(0)$,
 where each point in ${G_j}^{-1}(0)$ is counted multiplicity. We say that $i\sim j$ if $G_i=\gamma G_j$ for some $\gamma \in \mathcal P_c^{1}\setminus \{0\}.$ Therefore, we can split the set $\Lambda=\{1, \dots, q\}$ into disjoint equivalence classes $S_j$ such that $\Lambda=\cup_{j=1}^{l}S_j.$
First of all, assume that there exists $j\in \{1, \dots, l\}$ such that $S_j$ has at most $q-N-1$ elements. Put $R=\Lambda\setminus S_j,$ then $|R|\ge N+1.$ Let $s_0\in S_j$ and put $U=R\cup \{s_0\}.$ Without loss of generality, we may suppose that $U=\{s_0, s_1, \dots, s_{N+1}\}.$ Since hyperplanes $H_1, \dots, H_q$ are in $N$-subgeneral position in $\mathbb P^{n_{\mathcal D}}(\mathbb C),$ then exists $\gamma_0\in \mathbb C\setminus \{0\}$ and complex numbers $\gamma_j,j=1,\dots,N+1,$ are not simultaneous with zero such that $\sum_{j=0}^{N+1}\gamma_jH_{s_j}=0.$ Hence
$$ \sum_{j=0}^{N+1}\gamma_jH_{s_j}(\tilde{G})=\sum_{j=0}^{N+1}\gamma_jG_{s_j}\equiv 0.$$
By hypothesis of Theorem \ref{th3}, we see that all zeros of $\gamma_jG_j$ are forward invariant with respect to the translation $\tau(z)=z+c.$ We obtain a meromorphic mapping
$$ \mathcal G:=(\gamma_0G_{s_0}:\dots:\gamma_{N+1}G_{s_{N+1}}): \mathbb C^m\longrightarrow \mathbb P^{N+1}(\mathbb C)$$
with hyper-order $\varsigma (\mathcal G)<1.$ By Lemma \ref{lm23b},  we have $\gamma_0 G_{s_0}\equiv 0,$ then $H_{s_0}(\tilde{G})\equiv 0.$ Thus, the image $f(\mathbb C^m)$ is contained the hypersurface $D_{s_0}$ of $\mathbb P^{n}(\mathbb C)$.

Secondly, assume that $S_j$ has at least $q-N$ elements for all $j=1, \dots, l.$ Then $l\le \dfrac{q}{q-N}.$ Since $\{H_j\}_{j=1, \dots, q}$ is in $N$-subgeneral position, we can choose a subset $V\subset \{1,\dots, q\}$ with $|V|=n_{\mathcal D}+1$ such that $\{H_j\}_{j\in V}$ is linearly independent. Put $V_j=V\cap S_j,$ then we have $V=\cup_{j=1}^{l}V_j.$ Since each $V_j$ raise to $|V_j|-1$ equations over the field $\mathcal P_{c}^{1},$ then there are at least
$$ \sum_{j=1}^{l}(|V_j|-1)=|V|-l=n_{\mathcal D}+1-l\ge n_{\mathcal D}+1-\dfrac{q}{q-N}=n_{\mathcal D}-\dfrac{N}{q-N}$$
linearly independent relations over field $\mathcal P_{c}^{1}.$ Hence image of $G$ is contained the projective linear subspace over $\mathcal P_{c}^{1}$ with dimension at most $\Big[\dfrac{N}{q-N}\Big].$ Specially, if $q\ge 2N+1,$ then $G(z+c)\equiv G(z),$ this implies $f(z)\equiv f(z+c).$

\end{document}